\documentclass[12pt]{amsart}

\usepackage{amsmath}
\usepackage{latexsym}
\usepackage{amssymb}
\usepackage{amsfonts}
\usepackage{color,enumerate,tikz,tcolorbox}

\usetikzlibrary{decorations.pathreplacing,calc}

\usepackage{graphics}

\renewcommand{\proof}{\par\noindent{\it Proof.\ \ }}
\def\qed{\ifmmode\square\else\nolinebreak\hfill
$\Box$\fi\par\vskip12pt}

\def\l{\langle} \def\r{\rangle}

 \def\ZZ{{\mathbb Z}}

\def\calM{{\mathcal M}}

\def\ZZ{{\mathbb Z}}

\def\K{{\bf K}}

\def\Aut{{\rm Aut}}

\def\Cay{{\rm Cay}}
\def\Cos{{\rm Cos}}

\def\Circ{{\rm Circ}}

\def\D{{\rm D}}

\def\C{{\bf C}}

\def\Ga{{\it \Gamma}}
\def\Ome{{\it \Omega}}

\def\a{\alpha} \def\b{\beta}

\def\bfI{\mathbf{I}}
\def\la{\langle}
\def\ra{\rangle}
\def\lcm{{\rm lcm}}

\def\bI{\mathbf{I}}

\def\Sym{{\rm Sym}}

  \def\D{{\rm D}}

\def\Bicos{{\rm BiCos}}

\newtheorem{theorem}{Theorem}[section]%
\newtheorem{lemma}[theorem]{Lemma}%
\newtheorem{corollary}[theorem]{Corollary}%
\newtheorem{proposition}[theorem]{Proposition}%
\newtheorem{definition}[theorem]{Definition}%
\newtheorem{problem}[theorem]{Problem}%
\newtheorem{example}[theorem]{Example}%
\newtheorem{construction}[theorem]{Construction}%

\newtheorem{remark}[theorem]{Remark}

\begin{document}

\title{The graphs with a symmetrical Euler cycle}
\thanks{2020 MR Subject Classification: 20B25, 05C25, 05C35.}

\thanks{This work was partially supported by Australian Research Council Discovery 
project~DP160102323, and NSFC projects nos.~11771200, 11931005, 61771019, and NSFS no.~ZR2020MA044.
}

\thanks{Corresponding author: Cheryl E. Praeger}


\author{Jiyong Chen}
\address{School of Mathematical Sciences\\Xiamen University\\
Xiamen\\P. R. China }
\email{chenjy1988@xmu.edu.cn}

\author{Cai Heng Li}
\address{Department of Mathematics\\
Southern University of Science and Technology\\
Shenzhen\\P. R. China }
\email{lich@sustech.edu.cn}

\author{Cheryl E. Praeger}
\address{Department of Mathematics and
Statistics\\ The University of Western Australia\\ Crawley 6009, WA\\
Australia }
\email{cheryl.praeger@uwa.edu.au}

\author{Shu-Jiao Song}
\address{School of mathematics and information science\\ 
Yantai University\\
Yantai\\ P. R. China}
\email{shujiao.song@hotmail.com}

\date\today

\begin{abstract}

\ \\
\begin{center}
Dedicated to our friend and colleague Marston Conder on the occasion of his 65th birthday.
\end{center}

The edges surrounding a face of a map $\calM$ form a cycle $C$, called the boundary cycle of the face, and $C$ is often not a simple cycle.
If the map $\calM$ is arc-transitive, then there is a cyclic subgroup of automorphisms of $\calM$ which leaves $C$ invariant and is bi-regular on the edges of the induced subgraph $[C]$; that is to say, $C$ is a symmetrical Euler cycle of $[C]$.
In this paper we  determine the family of graphs (which may have multiple edges) whose edge-set can be sequenced to form a symmetrical Euler cycle. 
We first classify all graphs  and which have a cyclic subgroup of automorphisms acting bi-regularly on edges.
We then apply this classification to obtain the graphs possessing a symmetrical Euler cycle, and therefore are the (only) candidates for the induced subgraph of the boundary cycle 
of a face in an arc-transitive map.
\end{abstract}

\maketitle

\section{Introduction}

The graphs studied in this paper are finite, undirected, without loops, but may have multiple edges. Thus a graph $\Gamma=(V, E, \bI)$ consists of finite sets $V$ of vertices, and $E$ of edges, together with an incidence relation $\bI\subseteq V\times E$ such that each edge $e$ is incident with exactly two distinct vertices. We often suppress $\bI$ in the notation are write simply $\Gamma=(V, E)$. An edge $e$ of $\Ga$ incident with the two vertices $\a$ and $\b$ is sometimes denoted by $[\a,e,\b]$.
 Many graphs of this type admit natural embeddings as maps into closed surfaces: perhaps the simplest being the graph $\K_2^{(\lambda)}$ with exactly two vertices $V=\{\a,\b\}$ and $\lambda$ edges forming $E=\{[\a,e_i,\b]\mid 1\leq i\leq \lambda\}$, and embedded into a sphere with $\a,\b$ at the poles, and the edges $e_i$ arranged as $\lambda$ lines of `longitude' joining the two poles.  
 
 Motivating our investigation was our wish to study  embeddings of graphs from this family as arc-transitive maps in surfaces. If such a map has at least two faces then 
the edge-sequence $C$ obtained by moving around the boundary of a face forms a 'cycle', as defined in \eqref{e:cycle}. We show in Lemma~\ref{map-cycle} that this cycle is a `symmetrical Euler cycle' for the induced subgraph $[C]$: this is a cycle which admits a large subgroup of the corresponding dihedral group acting with at most two orbits on edges, see Subsection~\ref{s:intromaps} and Section~\ref{s:cycles} for more details.   
The most natural example for the cycle $C$ is the sequence obtained by traversing the edges around a simple $n$-cycle several times, say $\lambda$ times, and we call the induced subgraph $[C]$ in this example $\C_n^{(\lambda)}$, see Section~\ref{s:cycles}. Maps for which all boundary cycles are of this form with $\lambda=1$ (simple cycles) have been studied in  \cite{ci-map,C-map,s-map}. 
On the other hand, quite different subgraphs $[C]$ have been identified for maps with a single face in   \cite{1-face,Singerman}.  
The problem, which we address in the paper, is to determine the kinds of subgraphs $[C]$ that arise, induced by symmetrical Euler cycles $C$, and to describe the possible groups induced on these cycles by the subgroup of automorphisms leaving $[C]$ invariant.

Note that the subgraph induced by the edges of a boundary cycle of a map is connected, and for an arc-transitive map, the group induced on the cycle contains a cyclic subgroup having at most two edge-orbits which acts faithfully on each of these edge-orbits (Lemma~\ref{map-cycle}). 
In our first main result, Theorem~\ref{main-thm}, we broaden the scope of this study slightly, and  classify all connected graphs admitting  a cyclic group with at most two edge-orbits and acting faithfully on each of its edge-orbits. We find a dozen infinite families of examples. Then, in Theorem~\ref{thm-2}, we show that only six of these families contain graphs for which the edge set can be sequenced into a cycle preserved by a cyclic group with at most two  edge-orbits (of equal size).  To assist with our analysis we develop, in Section~\ref{s:coset}, the theory of coset graphs which may have multiple edges and which admit an edge-transitive group with two vertex-orbits. 

\subsection{Graphs admitting cyclic edge regular or bi-regular groups}

A graph $\Ga=(V,E)$ is a \emph{simple graph} if, for all distinct $\a,\b\in V$, the 
number of edges incident with both $\a$ and $\b$ is $0$ or $1$. 
Given a graph $\Ga=(V,E)$ and a positive integer $\lambda$, the {\it $\lambda$-extender} of $\Ga$ is the graph $\Ga^{(\lambda)}$ with vertex set $V$ such that each edge 
$[\a,e,\b]$ of $\Ga$ is replaced by $\lambda$ edges  $[\a,e_i,\b]$ ($1\leq i\leq 
\lambda$) of $\Ga^{(\lambda)}$. 
If $\Ga$ is simple, then $\Ga$ is said to be the {\it base graph} of $\Ga^{(\lambda)}$; in this case if $\a, \b$ are adjacent in $\Ga$, that is, if there exists an edge $[\a,e,\b]$ in $\Ga$, then there are exactly $\lambda$ edges of $\Ga^{(\lambda)}$ incident with $\a$ and $\b$, and we say that $\Ga$ has \emph{edge-multiplicity} $\lambda$.  For each edge $[\a,e,\b]$, we have 
$[\a,e,\b]=[\b,e,\a]$ and the edge $e$ corresponds to two arcs 
$(\a,e,\b)$ and $(\b,e,\a)$.

Let a group $G$ act on a set $\Ome$.
Then $G$ is called {\it transitive} or {\it bi-transitive} if $G$ has a single orbit or exactly two orbits in $\Ome$, respectively.
Further, if $G$ is finite, then $G$ is said to be {\it regular} or {\it bi-regular} on $\Omega$ if the permutation group $G^\Omega$ induced by $G$ on $\Omega$ is transitive and $|G^\Omega|=|\Ome|$, or $G^\Omega$ is bi-transitive and $|G^\Omega|=\frac{1}{2}|\Ome|$, respectively; in other words, $G^\Omega$ has at most two orbits in $\Omega$ and is faithful and regular on each.

An automorphism of a graph $\Ga=(V,E,\bI)$ is a permutation of $V\cup E$ which preserves $V$, $E$ and the incidence relation $\bI$. The set of automorphisms forms the automorphism group $\Aut\Ga$. Usually $\Aut\Gamma$ acts faithfully on $E$, see Lemma~\ref{l:aut} for details: the unique exceptions among connected graphs are the graphs $\K_2^{(\lambda)}$ mentioned above. Our first result Theorem~\ref{main-thm} presents a classification of connected graphs admitting a cyclic subgroup of automorphisms that is regular or bi-regular on edges. The exceptional graphs $\K_2^{(\lambda)}$ are treated separately in detail in Proposition~\ref{p:k2}. The families of examples, apart from $\K_2^{(\lambda)}$, are defined in Section~\ref{sec:examples}. Note that for graphs $\Gamma=(V,E)$ with an edge partition $E =E_1\cup E_2$, we sometimes write 
$\Gamma = [E_1] + [E_2]$ to give a rough description of the graph structure, even though this notation does not uniquely define the graph in general, see Subsection~\ref{r:sum}. We give a precise description  in Section~\ref{sec:examples} of all the graphs in the tables for Theorem~\ref{main-thm}.

\begin{theorem}\label{main-thm}
Let $\Ga=(V,E)$  be a connected graph with $|V|\geq3$ such that a cyclic subgroup $G\leq \Aut\Ga$ is regular or bi-regular on $E$, and has $N_V$ orbits on $V$. Then
$\Gamma=\Gamma_0^{(\lambda)}$ and $|G|=\lambda N$, for some $\lambda, N$, and either 

 \begin{enumerate}[{\rm (I)}]
  	\item $G$ is regular on $E$ and $\Gamma_0, N, N_V$ are as in one of the lines of Table~\ref{t:mainreg}, or 
  	
  	\item $G$ is bi-regular on $E$ and  $\Gamma_0, N$ are as in one of the lines of 
  	Table~\ref{t:mainbireg1} if $N_V=1$, and $\Gamma_0, N, N_V$ are as in one of the lines of Table~\ref{t:mainbireg2} if $N_V\geq2$. In particular $N_V\leq 3$.  
  \end{enumerate}
\end{theorem}
 
\begin{center}
\begin{table}
\begin{tabular}{cccll}
$\Gamma_0$ & $N$ & $N_V$ & Conditions & Reference  \\ 
\hline 
$\C_n$ 		& $n$ 	& $1$ &  $n\geq 3$		& Lemma~\ref{l:n-cycle}  \\ 
$\K_{s,t}$ 	& $st$ 	& $2$ & $\gcd(s,t)=1$, $st>1$ & Definition~\ref{d:kst}, Lemma~\ref{l:kst} \\ 
\hline 
\end{tabular} 
\caption{Table for Theorem~\ref{main-thm} with $G$ regular on $E$}\label{t:mainreg}
\end{table}
\end{center}

\begin{center}
\begin{table}
\begin{tabular}{ccll}
$\Gamma_0$ & $N$  & Conditions & Reference  \\ 
\hline 
$\C_{n}^{(2)}$& $n$ 	&  $n\geq 3$& Lemma~\ref{l:n-cycle}  \\ 
$\C_{2n}+n\K_2^{(2)}$& $2n$ 	&  $n\geq 2$& Corollary~\ref{c:circ1}(a)  \\ 
$2\C_{n}+n\K_2^{(2)}$& $2n$ 	&  $n\geq 3$ odd& Corollary~\ref{c:circ1}(b)  \\ 
$\Circ(n,S)$& $n$ 	&  $n\geq 5$, $S=\{a,-a,b,-b\}$,& Definition~\ref{d:circ1}  \\ 
		& 	 	& $|S|=4$, $\gcd(n,a,b)=1$&   \\ 
\hline 
\end{tabular} 
\caption{Table for Theorem~\ref{main-thm} with $G$ bi-regular on $E$ and transitive on $V$}\label{t:mainbireg1}
\end{table}
\end{center}

\begin{center}
\begin{table}
\begin{tabular}{lccll}
$\Gamma_0$ & $N$  & $N_V$ & Conditions & Reference  \\ 
\hline 
$\C_{n}=$	& $n$ &$2$	&  $n\geq 4$ even& Lemma~\ref{l:n-cycle}  \\ 
\ \ $\frac{n}{2}\K_2+\frac{n}{2}\K_2$ &&& \\
$\K_{s,t}^{(2)}$ 	& $st$ 	 &$2$	& $\gcd(s,t)=1$, $st>1$ & Lemma~\ref{l:kst} \\
$\C_{2r}[s\K_1,t\K_1]$& $rst$  &$2$	&  $r\geq 2$, $\gcd(s,t)=1$, $st>1$& Lemma~\ref{l:cyclest}  \\ 
$r\K_{2}^{(2t)}+\K_{2r,t}$& $2rt$  &$2$		&  $rt\geq1$, $\gcd(2r,t)=1$& Lemma~\ref{l:ck}  \\ 
$r\K_{2}^{(2t)}+2\K_{r,t}$& $2rt$  &$2$		&  $r$ odd, $rt\geq1$, $\gcd(r,t)=1$& Lemma~\ref{l:ck2}  \\ 
$r\C_{n}^{(t)}+\K_{nr,t}$& $nrt$  &$2$		&  $n\geq 3$, $rt\geq1$, $\gcd(nr,t)=1$& Lemma~\ref{l:ck}  \\ 
$r\C_{su}^{(t)}+u\K_{sr,t}$& $srut$  &$2$		& $su\geq3$,  $u\geq2$,  $\gcd(r,u)=1$& Lemma~\ref{l:ck2}  \\ 
					& &		& and $\gcd(sr,t)=1$&   \\ 
\hline
$r\K_{sr',ut'}^{(t)} + $ 	& $rr'stt'u$ 	 &$3$	& $\gcd(r,r')=\gcd(t,t')=1$ & Lemma~\ref{l:kk} \\
\ \ $ r'\K_{sr,ut}^{(t')}$ 	& &	&  and $\gcd(sr,ut)=\gcd(sr',ut')=1$ &  \\
\hline 
\end{tabular} 
\caption{Table for Theorem~\ref{main-thm} with $G$ bi-regular on $E$ and intransitive on $V$}\label{t:mainbireg2}
\end{table}
\end{center}

\subsection{Cycles in graphs and boundary cycles of maps}\label{s:intromaps}

A {\it cycle of length $\ell$} in a graph $\Ga=(V, E)$, sometimes called an \emph{$\ell$-cycle}, is a sequence
\begin{equation}\label{e:cycle}
C=(e_1,e_2,\ldots,e_{\ell})
\end{equation}
of $\ell$ pairwise distinct edges each of the form
$[\a_{i-1},e_i,\a_i]$, for $1\leq i\leq \ell$, 
and we read the subscripts modulo $\ell$ so that, in particular,  $\alpha_\ell=\alpha_0$.
The {\it edge induced subgraph} $[C]$ of a cycle $C$ is the graph with vertex set $V(C):=\{\alpha_i\mid 1\leqslant i\leqslant \ell\}$, edge set $E(C):=\{e_i\mid 1\leqslant i\leqslant \ell\}$, and incidence as in $\Ga$.
We call $C$ an {\it Euler cycle} of $\Ga$ if $E(C)=E$, that is, the cycle `passes through each edge of $\Ga$ exactly once'. 
If $\Ga$ possesses an Euler cycle then, in particular, $\Ga$ is connected. 
Further, if $C$ as in \eqref{e:cycle} is an Euler cycle of $\Ga$, then since each edge of $E$ occurs exactly once in $C$, and since whenever a vertex $\a=\a_i$ then both $e_i$ and $e_{i+1}$ are incident with $\a$, it follows that \emph{for each vertex $\a\in V$, the number of edges (of $C$, and hence of $\Gamma$) incident with $\a$ is even.} 

For any cycle $C$, the following bijections on $E(C)$ form a dihedral group of 
order $2\ell$:
\begin{equation}\label{e:varphi}
\varphi:e_i\to e_{i+1} \quad\mbox{and}\quad \tau: e_i\to e_{\ell+1-i}\quad \mbox{(for each $i$,  reading subscripts modulo $\ell$)}.
\end{equation}
The reflection $C\tau=(e_{\ell},e_{\ell-1},\ldots,e_1)$, and each shift $C\varphi^i=(e_{i+1},\ldots,e_{\ell},e_1,\ldots,e_{i})$ are also $\ell$-cycles of $\Ga$, and
\begin{equation}\label{e:DC}
\D(C):= \l \varphi,\tau\r=\la \varphi,\tau\mid\varphi^\ell=\tau^2=1, \varphi^\tau=\varphi^{-1}\ra \cong \D_{2\ell}
\end{equation}
is the subgroup of all permutations of $E(C)$ which preserve the edge-sequencing of $C$, up to rotations and reflections. Thus the subgroup of $\Aut\Ga$ leaving the cycle $C$ invariant induces on $[C]$ a subgroup $H(C)$ of $\Aut[C]$ contained in $\D(C)$. If this subgroup contains $\varphi^2$, that is to say, if there is an element $\sigma \in\Aut(\Ga)$ such that    
$\sigma|_{[C]}=\varphi^2$, then we say that $C$ is {\it symmetrical in $\Ga$}. 
In particular, for an Euler cycle $C$ of a graph $\Ga=(V,E)$, $C$ is symmetrical if and only if $\Aut\Ga$ contains a cyclic subgroup which preserves $C$ and is regular or bi-regular on $E$. For future reference, we note the definitions of the elements $\varphi^2$ and $\varphi\tau$ of $D(C)$:
\begin{equation}\label{e:varphitau}
\varphi^2:e_i\to e_{i+2} \quad\mbox{and}\quad \varphi\tau: e_i\to e_{\ell-i}\quad \mbox{(for each $i$,  reading subscripts modulo $\ell$)}.
\end{equation}

For the exceptional graphs $\K_2^{(\lambda)}$, we show in Proposition~\ref{p:k2} that there is a symmetrical Euler cycle if and only if $\lambda$ is even. 
For all other (connected) graphs we apply Theorem~\ref{main-thm} to determine whether or not they have a symmetrical Euler cycle.  

\begin{theorem}\label{thm-2}
Let $\Ga=(V,E)$ be a graph with $|V|\geq3$ which has a symmetrical Euler cycle $C$, and suppose that  the subgroup of $\Aut\Ga$ which preserves $C$ induces the subgroup $H(C)$  of $D(C)$, as in \eqref{e:DC}.
\begin{enumerate}
\item[(a)] Then $(\Gamma, H(C))$ are as in one of the lines of Table~\ref{t:main-cycle}.

\item[(b)] Conversely, if $\Ga$ is one of the graphs in lines $1-5$ of Table~\ref{t:main-cycle}, then $\Gamma$ has a symmetrical Euler cycle $C$ with the given $H(C)$; and also if $\Gamma$ is as in line $6$ of Table~\ref{t:main-cycle} with at least one of $\gcd(n,a+b)=1$ or $\gcd(n,a-b)=1$, then $\Gamma$ has a symmetrical Euler cycle.
\end{enumerate}
\end{theorem}

\begin{center}
\begin{table}
\begin{tabular}{clll}
$\Gamma$ & Conditions for $C$ to exist & $H(C)$ & Reference  \\ 
\hline 
$\C_{n}^{(\lambda)}$& all $n\geq3, \lambda\geq1$ &  $H(C)=D(C)$& Lemma~\ref{l:kc-euler}(a) \\ 
$\K_{s,t}^{(\lambda)}$&  $st>1$, $\gcd(s,t)=1$,  $\lambda$ even &  $H(C)=\l\varphi^2,\varphi\tau\r$& Lemma~\ref{l:kc-euler}(b)  \\ 
$(\C_{2n}+n\K_2^{(2)})^{(\lambda)}$& $n\geq2$, $n$ even, $\lambda\geq1$& $H(C)=\l\varphi^2,\varphi\tau\r$ & Lemma~\ref{Cyc-2-E}(a)  \\ 
$(2\C_{n}+n\K_2^{(2)})^{(\lambda)}$& $n\geq3$, $n$ odd, $\lambda\geq1$& $H(C)=\l\varphi^2,\varphi\tau\r$ & Lemma~\ref{Cyc-2-E}(b)  \\ 
 $(\C_{2r}[s\K_1,t\K_1])^{(\lambda)}$ & $r\geq2, st>1$, $\gcd(s,t)=1$, $\lambda\geq1$	& $H(C)=\l\varphi^2,\tau\r$&  Lemma~\ref{l:kst-euler}\\
$(\Circ(n,S))^{(\lambda)}$&  $n\geq 5$, $S=\{a,-a,b,-b\}$, & $H(C)\geq\l\varphi^2\r$ 	&  Lemma~\ref{Cyc-2-E}(c)   \\ 
		& $|S|=4$,  $\gcd(n,a,b)=1$,  $\lambda\geq1$	& &   \\ 
\hline 
\end{tabular} 
\caption{Table for Theorem~\ref{thm-2} with $C$ a symmetrical Euler cycle for $\Gamma$}\label{t:main-cycle}
\end{table}
\end{center}

As we discuss in Remark~\ref{r:circ-rem}, if $\Gamma=\Gamma_0^{(\lambda)}$ with $\Gamma_0=\Circ(n,S)$ as in line 6 of Table~\ref{t:main-cycle}, and if $\Gamma_0$ is a \emph{normal Cayley graph} (that is the translation subgroup $\ZZ_n$ is normal in $\Aut\Gamma_0$), then the  the condition `$\gcd(n,a+b)=1$ or $\gcd(n,a-b)=1$' is necessary and sufficient for existence of a symmetrical Euler cycle $C$. A complete analysis of these graphs would need a better understanding of any non-normal, edge-transitive, Cayley graphs in this family.

The principle motivation for our work was the study of boundary cycles of faces of arc-transitive maps. The link between arc-transitive maps and symmetrical Euler cycles is made explicit in Lemma~\ref{map-cycle}. A natural problem which we plan to explore in further work is to understand which of these graphs from Theorem~\ref{thm-2} actually arise in arc-transitive maps.

\begin{problem}\label{pr1}
Determine which of the graphs in Table~\ref{t:main-cycle} arise as the induced subgraph of the boundary cycle of a face in an arc-transitive map. 
\end{problem}

We were able to find a partial answer in our study of vertex-rotary maps, that is, arc-transitive maps for which a vertex-stabiliser is cyclic and regular on the edges incident with it.  We prove in \cite[Theorem 1.7]{LPS-rotary} that each face boundary cycle $C$ for a vertex-rotary map has induced subgraph $[C]=\C_n^{(\lambda)}$ as in line ~1 of Table~\ref{t:main-cycle}, and moreover the cycle $C$ is of the form 
$C_0^{(\lambda)}$ (as in Proposition~\ref{p:biregext}) for a simple $n$-cycle $C_0$. 
In addition, for $n\equiv 2\pmod{4}$, we construct in \cite[Section 6]{LPS-rotary} infinitely many examples of such maps from extenders of complete bipartite graphs, namely, if $n=2m$ with $m$ odd,  and if $\lambda$ is such that $\gcd(n,\lambda)=2$, we construct an embedding of $\K_{m,m}^{(2\lambda)}$ with all face boundary cycles of the form $C_0^{(\lambda)}$ for a simple $n$-cycle $C_0$. In particular, it would be interesting to know precisely which values of $n, \lambda$ are possible for vertex-rotary maps.    


\section{Preliminaries and the exceptional graph $\K_2^{(\lambda)}$}\label{s:prelim}

\subsection{Notation}\label{r:sum}

For a graph $\Ga=(V,E,\bI)$, and any subset of edges $E'\subseteq E$, the \emph{edge-induced subgraph} is the graph $[E']:=(V(E'), E', \bfI')$, where 
\begin{center}
$V(E')=\{ \a\mid (\a,e)\in \bI, \mbox{ for some }e\in E'\}$, and $\bfI'= \bfI\cap (V(E')\times E'$.
\end{center}
An \emph{isolated vertex} of a graph $\Ga$, is a vertex which is not incident to any edge  of $\Gamma$. Thus $\Gamma$ has no isolated vertices if and only if if $V(E)=V$. The graph $\Gamma$ is \emph{connected} if, for each pair of vertices $\a,\b$ there exists a sequence of edges $e_1,\dots, e_r$ such that $e_i=[\a_{i-1},e_i,\a_i]$ for each $i$, and $\a_0=\a$ and $\a_r=\b$.  
A \emph{connected component of $\Gamma$} is either (i) a one-vertex graph $(\{\a\},\emptyset)$ where $\a$ is an isolated vertex, or (ii) a connected edge-induced subgraph 
$[E']$ for some non-empty $E'\subseteq E$ such that, for each edge $e \in E$,  either  $e\in E'$, or $e$ is incident with no vertex of $V(E')$. 

For some of the graphs in Theorem~\ref{main-thm}, the group considered has two edge-orbits, and we view the graph as an edge-disjoint union of two smaller graphs. We use the following notation: given graphs $\Gamma_i=(V_i, E_i, \bI_i)$, for $i=1,2$, where the vertex sets $V_1, V_2$ may overlap but $E_1\cap E_2=\emptyset$, the \emph{edge-disjoint union} of $\Gamma_1$ and $\Gamma_2$ is the graph $\Gamma_1 + \Gamma_2 =(V, E, \bI)$, such that
\begin{center}
$V=V_1\cup V_2$,\ $E=E_1\cup E_2$,\ and\ $\bI=\bI_1\cup\bI_2\subset V\times E$. 
\end{center}

\subsection{The graph $\K_2^{(\lambda)}$}

First we prove the assertion, mentioned in the introduction, that essentially the only graphs for which the automorphism group is not faithful on edges are those with at least one component $\K_2^{(\lambda)}$.

\begin{lemma}\label{l:aut}
Let $\Gamma=(V, E, \bfI)$ be a graph with no isolated vertices and $G\leq \Aut\Gamma$. Then either $G$ acts faithfully on $E$, or some connected component of $\Gamma$ is isomorphic to $K_2^{(\lambda)}$, for some $\lambda$. 
\end{lemma}

\proof
Let $K=G_{(E)}$ be the kernel of the action of $G$ on $E$. By the definition of an automorphism, $K$ acts faithfully on $V$. If $K$ is trivial on $V$ there is nothing to prove, so assume that $x\in K$ and $\a\in V$ such that $\a^x\ne\a$. Since $\Gamma$ has no isolated vertices there exists an edge $e=[\a,e,\b]$, and since $x$ fixes $e$ it follows that $x$ interchanges $\a$ and $\b$. Since $K$ fixes all edges, it follows that each edge incident with $\a$ is also incident with $\b$, and conversely, and hence $\Gamma$ has a connected component with vertex set $\{\a,\b\}$, and it is isomorphic to 
$K_2^{(\lambda)}$, for some $\lambda$.  
\qed

The graph $\Gamma=\K_2^{(\lambda)}=(V,E)$ has automorphism group $\Aut\Gamma= \Sym(E)\times\Sym(V)\cong \Sym(\lambda)\times \ZZ_2$, with $\Sym(\lambda)$ acting naturally on $E$ and fixing $V$ pointwise, and $\Sym(V)$ acting naturally on $V$ and fixing $E$ pointwise. The additional automorphisms give extra examples of edge-regular and bi-regular actions, and it is instructive to consider these graphs separately here.

\begin{proposition}\label{p:k2}
Let $\Gamma=\K_2^{(\lambda)}=(V,E)$ for some $\lambda\geq1$, and let $G=\l g\r\leq \Aut\Gamma$ be  regular or bi-regular on $E$. Then  
\begin{enumerate}
\item[(a)] $G\leq L\times\l y\r$, where $\l y\r=\Sym(V)\cong \ZZ_2$, $L=\l x\r\cong\ZZ_\lambda$ is a cyclic edge-regular subgroup of $\Sym(E)$, and  $G, \lambda$ are as in one of the lines of Table~\ref{t:k2}.

\item[(b)] Moreover, $\Gamma$ has a symmetrical Euler cycle $C$ if and only if $\lambda$ is even, and in this case $H(C)=\l\varphi, \tau\r=D(C)$, with $\varphi,\tau, D(C)$ as in \eqref{e:varphi} and \eqref{e:DC}.

\end{enumerate}
\end{proposition}

\begin{center}
\begin{table}
\begin{tabular}{ccccl}
$G$ & $|G|$ & Action on $V$ & Action on $E$ & Conditions/comments  \\  \hline 
$\l x\r\times \l y\r$& $2\lambda$ 	& transitive & regular &  $\lambda$ odd, $G_{(E)}\ne 1$ \\ 
$\l xy\r$       & $\lambda$ 	& transitive & regular &  $\lambda$ even \\ 
$\l x\r$        & $\lambda$ 	& trivial    & regular &  $\lambda$ arbitrary \\ 
$\l x^2\r\times \l y\r$& $\lambda$ 	& transitive & bi-regular &  $\lambda=2\lambda_0$, $\lambda_0$ odd, $G_{(E)}\ne 1$ \\ 
$\l x^2y\r$       & $\lambda_0$ 	& transitive & bi-regular &  $\lambda=2\lambda_0$, $\lambda_0$ even \\ 
$\l x^2\r$        & $\lambda_0$ 	& trivial    & bi-regular &  $\lambda=2\lambda_0$, $\lambda_0$ arbitrary \\ 
\hline 
\end{tabular} 
\caption{Edge regular and bi-regular actions on $\K_2^{(\lambda)}$}\label{t:k2}
\end{table}
\end{center}

\proof
(a) The projection $\pi:G\to \Sym(E)$ has image a cyclic edge-regular or bi-regular subgroup, and in either case there exists a cyclic edge-regular subgroup $L\leq\Sym(E)$ containing $\pi(G)$. Thus $G\leq L\times\l y\r$ with $L=\l x\r\cong\ZZ_\lambda$. 
Suppose first that $G$ is edge-regular, so $\pi(G)=L$ and $|G|$ is a multiple of $\lambda$. One possibility is that $G=L\times \l y\r$, and since $G$ is cyclic this implies that $\lambda$ is odd and line 1 of Table~\ref{t:k2} holds. So suppose that $G$ is a proper subgroup of $L\times \l y\r$. Since $|G|$ is a multiple of $\lambda$, it follows that $|G|=\lambda$.  If $G$ is vertex-transitive then $G=\l xy\r\cong \ZZ_\lambda$ and this implies that $\lambda$ is even, and line 2 of Table~\ref{t:k2} holds. On the other hand, if $G$ is intransitive, then $G=\l x\r\cong \ZZ_\lambda$ and line 3 of Table~\ref{t:k2} holds (and here $\lambda$ is arbitrary).

Suppose now that $G$ is edge-bi-regular, so $\lambda=2\lambda_0$ is even and $\pi(G)=\l x^2\r\cong\ZZ_{\lambda_0}$, and $G\leq \l x^2\r\times \l y\r$. If $G= \l x^2\r\times \l y\r$, then $\lambda_0=|x^2|$ is odd since $G$ is cyclic, and line 4 of Table~\ref{t:k2} holds. Suppose now that $G$ is a proper subgroup of  $\l x^2\r\times \l y\r$, so $|G|=\lambda_0$. If $G$ is vertex-transitive then $G=\l x^2y\r\cong \ZZ_{\lambda_0}$ and this implies that $\lambda_0$ is even, and line 5 of Table~\ref{t:k2} holds. On the other hand, if $G$ is intransitive, then $G=\l x^2\r\cong \ZZ_{\lambda_0}$ and line 6 of Table~\ref{t:k2} holds (and here $\lambda_0$ is arbitrary).

(b) By our comments above in an Euler cycle, each vertex is incident with an even numberof edges, so there is no Euler cycle if $\lambda$ is odd. Suppose that $\lambda$ is even. Then it is easy to sequence the edges into a symmetrical Euler cycle using, say, the edge-regular automorphism $x$: choose $e\in E$ and define $e_0=e$ and $e_i=e^{x^i}$ for $1\leq i\leq \lambda-1$. Then $C=(e_0,\dots,e_{\lambda-1})$ is an Euler cycle and $x$ induces $\varphi\in D(C)$ and fixes each of the vertices. Also there is an involution $z\in N_{\Sym(E)}(\l x\r)$ given by $z:e_i\to e_{\lambda-i}$ for each $i$, and hence $z$ induces $\tau$ in $H(C)$ so $H(C)=D(C)$.
\qed

\section{Cycles, symmetrical Euler cycles, and maps}\label{s:cycles}

Let $C$ be an $\ell$-cycle in a graph $\Gamma=(V, E, \bI)$, as in \eqref{e:cycle}, with edges $[\a_{i-1},e_i,\a_i]$, reading the subscripts modulo $\ell$. If the vertices 
$\a_i$ are pairwise distinct, then $[C]$
is a simple graph of valency $2$, and hence $[C]=\C_\ell$. In this case we call $C$  a \emph{simple $\ell$-cycle}, and since the vertices $\a_i$ are pairwise distinct, 
we can identify $C$ with the vertex sequence $(\a_1,\a_2,\dots,\a_{\ell})$.

Other examples of $\ell$-cycles $C$ arise naturally with $[C]=\C_n^{(\lambda)}$ where $\ell=n\lambda$, with the edges sequenced in any way so that $\alpha_i=\alpha_j$ whenever $i\equiv j\pmod{n}$.
The set of all $\ell$-cycles with this property forms a single orbit under the automorphism group of $\C_n^{(\lambda)}$, and we
call each such $\ell$-cycle \emph{standard}, or a \emph{standard edge-sequencing}.
However not all cycles in graphs are of this kind. Even for  $[C]=\C_n^{(\lambda)}$ with $\lambda>1$, there are many non-standard ways to
seqence the edges to form a cycle.
For example, if $e_{ij}^1, e_{ij}^2$ denote the two edges of $\C_3^{(2)}$ incident with vertices $i,j$, then
$(e_{12}^1, e_{12}^2, e_{23}^1, e_{13}^1, e_{13}^2, e_{23}^2)$ is a non-standard $6$-cycle.
Moreover the induced subgraph $[C]$ of a cycle $C$ may be quite different from $\C_n^{(\lambda)}$, as the next examples show (some of which appear again in our main theorems). 

\begin{example}\label{Ex-1}
{\rm Let $\Ga=\Sigma^{(\lambda)}$, with $\lambda=2$ and with $\Sigma$ 
one of the following simple graphs.
Then $\Ga$ has constant edge-multiplicity $2$, and its edge set can be sequenced to form an Euler cycle. (Finding these sequencings is left as an easy exercise.)

\begin{enumerate}

\item[(a)] $\Sigma$ is the complete bipartite graph $\K_{2,3}$, so $\Ga$ has $\ell=12$ edges and $n=5$ vertices;

\item[(b)] $\Sigma$ is the cartesian product ${\C_{r}}\,\Box\,\K_{2}$, so $\Ga$ has $\ell=6r$ edges and $n=2r$ vertices, and $\gamma$ is regular of valency $3$;

\item[(c)] $\Sigma = (V,E,I)$ with $V=\mathbb{Z}_r\times \mathbb{Z}_2$ and $E = \{ e_i, f_i\mid i\in\mathbb{Z}_r\}$ such that, for each $i\in\mathbb{Z}_r$, $e_i$ is incident with $(i,0)$ and $(i+1,0)$, and $f_i$ is incident with $(i,0)$ and $(i,1)$. Here   $\Ga=\Sigma^{(2)}$ has $\ell=4r$ edges, $n=2r$ vertices, and we note that $\ell = n\lambda$.

\end{enumerate}
}
\end{example}

These examples illustrate that, for an $\ell$-cycle $C$,  the induced graph $[C]$ need not be regular (examples (a) and (c)), or $[C]$ may be regular of valency greater than $2$ (example (b)), and even if $[C]$ has constant edge-multiplicity $\lambda$ and $n\lambda$ edges, where $n$ is the number of vertices, $[C]$ need not be $\C_n^{(\lambda)}$ (example (c) again). These observations hold in particular for the Euler cycles of the graphs in Example~\ref{Ex-1}.

\subsection{Symmetrical cycles and arc-transitive maps}

Let $C=(e_1,e_2,\ldots,e_\ell)$ be an $\ell$-cycle in a graph $\Gamma$. As we noted in Section~\ref{s:intromaps}, each rotation (shift) $C\varphi^i$ and reflection (reversal followed by a shift) $C\tau\varphi^i$ of $C$ is again a cycle involving the same edge-set $E(C)$ as $C$. This collection of cycles if called the \emph{sequence class} of cycles containing $C$. We say that an automorphism $g\in\Aut\Gamma$ \emph{preserves $C$} if $g$ leaves invariant the sequence class of $C$. Thus the subgroup $(\Aut\Gamma)_C$ of $\Aut\Gamma$ preserving $C$ induces a subgroup $H(C)$ of automorphisms of the induced subgraph $[C]$ such that $H(C)\leq \D(C)$ with $\D(C)$ as in \eqref{e:DC}. 
As discussed above, $C$ is called symmetrical if $H(C)$ contains $\l\varphi^2\r$, and for such cycles, Theorem~\ref{thm-2} determines all possible induced subgraphs $[C]$.
Apart from $\l\varphi^2\r$ itself, the subgroups $H(C)$  of $\D(C)$ containing $\l\varphi^2\r$ are listed in Table~\ref{t:map}. These subgroups are relevant for studying edge-transitive or face-transitive embeddings of graphs in Riemann surfaces. 

A {\it map} $\calM=(V,E,F)$ is a 2-cell embedding of a graph $\Ga=(V,E)$ in a closed surface with face set $F$, and the graph $\Ga$ is called the {\it underlying graph} of $\calM$. An {\it arc} of a map, or a graph, is a an incident vertex--edge pair, and 
a map $\calM$ is called {\it arc-transitive} if its automorphism group $\Aut\calM$ is transitive on the arc-set of $\calM$. The next result describes the the kinds of groups $H(C)$ which may arise for boundary cycles $C$ of arc-transitive maps. In \cite{BigMaps} we discuss finite maps $\calM=(V,E,F)$ with arc-transitive automorphism groups, where the underlying graph $\Gamma=(V,E)$ may have multiple edges, as in this paper. The only facts about maps we use in the proof of Lemma~\ref{map-cycle} are that each edge in $E$ is incident with exactly two faces in $F$, together with the arc-transitivity of the map group $\Aut\calM \leq \Aut\Gamma$.

\begin{table}
\begin{center}
\caption{Groups $H(C)$ for a boundary cycle $C$ of an  arc-transitive map $\calM$}\label{t:map}
\begin{tabular}{|c|c|c|l|}
\hline 
$H(C)$ & orbits in $E(C)$ & orbits in $V(C)$ & Comments \\  \hline 
$\D(C)$ & $E(C)$ & $V(C)$ &  \\  
$\l\varphi\r$ & $E(C)$ & $V(C)$ & equals $\l\varphi^2\r$ if $\ell$ odd \\  
$\la\varphi^2,\varphi\tau\ra$ & two orbits & $V(C)$ & $\ell$ even \\ 
$\la\varphi^2,\tau\ra$ & $E(C)$ & two orbits if $\lambda=1$ & $\ell$ even \\ 
\hline 
\end{tabular} 
\end{center}
\end{table}

\begin{lemma}\label{map-cycle}
Let $\calM=(V,E,F)$ and $G\leqslant\Aut\calM$ such that $G$ is arc-transitive on $\calM$, and each edge in $E$ is incident with exactly two faces in $F$.
Let $f\in F$ be a face, and $C=(e_1,\dots, e_\ell)$ the boundary cycle of $f$.
Then the stabiliser $G_f$ induces a subgroup $H(C)$ of $\D(C)$ as in one of the lines of Table~$\ref{t:map}$, which gives also the $H(C)$-orbits on edges $E(C)$ and vertices $V(C)$ of $C$. In each case, $H(C)$ contains the cyclic subgroup $\l\varphi^2\r$ which is 
regular (if $\ell$ is odd) or bi-regular (if $\ell$ is even) on $E(C)$, and $\l\varphi^2\r$  has at most two orbits in $V(C)$. In particular $C$ is a symmetrical Euler cycle of $[C]$.
\end{lemma}

\proof
For $C=(e_1,\dots, e_\ell)$ we have $e_i=[\a_{i-1},e_i,\a_i]$ for each $i$, reading subscripts modulo $\ell$. Let $f_i$ be the second face incident with $e_i$, for each $i$. The setwise stabiliser  $G_f$ preserves $C$ and induces a subgroup $H(C)$ of $\D(C)$.  As $\calM$ is $G$-arc-transitive, there exist $g,h\in G  $ such that $(\a_0,e_1,\a_1)^h=(\a_1,e_1,\a_0)$ and $(\a_0,e_1,\a_1)^g=(\a_1,e_2,\a_2)$. Suppose first that $f^h=f$. Then also $f_1^h=f_1$. For each $i$ there exists $g_i\in\Aut\calM$ such that $e_1^{g_i}=e_i$ and so $g_i$ maps the pair $\{f,f_1\}$ to the pair of faces $\{f, f_i\}$ incident with $e_i$. In particular $f^{g_i^{-1}}$ is one of $f, f_1$ and so is fixed by $h$ and hence $f^{g_i^{-1}hg_i} = f^{g_i^{-1}g_i}=f$. It follows that $H(C)$ contains each edge-reflection, and hence contains $\varphi^2$ as well as the edge-reflection induced by $h$. Thus in this case $H(C)$ is either  $\la\varphi^2,\varphi\tau\ra$ with $\ell$ even and has two orbits in $E(C)$, or is $\D(C)$, and in either case $H(C)$ is transitive on $V(C)$. This we may assume that each element $h$ such that $(\a_0,e_1,\a_1)^h=(\a_1,e_1,\a_0)$ interchanges $f$ and $f_1$, and hence that $g_i^{-1}hg_i$ interchanges $f$ and $f_i$, for each $i$. 
 If $g$ leaves the face $f$ invariant then $H(C)$ contains $\varphi$ so $H(C)=\l\varphi\r$ (by our assumption on the elements $h$) and $H(C)$ is transitive on both $E(C)$ and $V(C)$. 
 
 So assume also that $g$ does not fix $f$ for any such $g$ mapping $(\a_0,e_1,\a_1)$ to $(\a_1,e_2,\a_2)$.    Then $f_1^g=f, f^g=f_2$, and hence $hg$ leaves $f$ invariant. Now $hg$ induces a reflection of $C$ in the vertex $\a_1$. Similar arguements with the edges $e_1, e_2$ replaced by $e_2, e_3$ yield one of the groups above, for $H(C)$ or that $H(C)$ contains also the reflection of $C$ in the vertex $\a_2$. The product of these two reflections is a generator of the group $\l\varphi^2\r$, and hence these two reflections generate the group $H(C)=\l \varphi^2, \tau \r$, which is regular on $E(C)$ and has two vertex-orbits when $\ell$ is even. (If $\ell$ is odd the group generated is $\D(C)$.) 
\qed

\begin{remark}
{\rm
We note that the subgroup $\l\varphi^2,\tau\r$ contains all reflections of $C$ `through a vertex': $\tau$ reflects $C$ through the vertex $\a_0$ (between the edges $e_\ell$ and $e_1$), while $\tau\varphi^2$ reflects $C$ through the vertex $\a_1$ (between the edges $e_1$ and $e_2$), etc. Thus, for a symmetrical cycle $C$ as in \eqref{e:cycle}, to prove that $H(C)=\l \varphi^2,\tau\r$ we need only prove that $H(C)$ contains one of these reflections. Similar comments apply to the subgroup $\l\varphi^2,\varphi\tau\r$, which contains all reflections of $C$ `through an edge', where $\varphi\tau$ reflects $C$ through the edge $e_\ell$, etc. To prove that $H(C)=\l \varphi^2,\varphi\tau\r$, for  a symmetrical cycle $C$, we need only prove that $H(C)$ contains one reflection through an edge. We note also that it is sometimes more convenient in our analysis to label the edges of an $\ell$-cycle as $(e_0, e_1,\dots, e_{\ell-1})$ for working with the action of the cyclic group $\ZZ_\ell$.
}
\end{remark}

\subsection{Symmetrical Euler cycles and $\lambda$-extenders}\label{s:symext}

It turns out that if a  graph $\Gamma$ (not necessarily a simple graph) has a symmetrical Euler cycle, then so does its $\lambda$-extender for each $\lambda$. 
We prove this as a consequence of a more general result about cyclic edge regular or bi-regular subgroups of automorphisms. It is useful to use the following notation: for  $\Ga=(V,E)$, we have $\Ga^{(\lambda)}=(V,E^{(\lambda)})$, and we label the $\lambda$ edges in $E^{(\lambda)}$ corresponding to the edge  $e=[\b,e,\a]\in E$ by
\begin{equation}\label{e:Elambda}
e^j=[\b,e^j,\a],\ \mbox{for $1\leqslant j\leqslant \lambda$.}
\end{equation}

\begin{proposition}\label{p:biregext}
Suppose that, for a graph $\Ga=(V,E)$, $\Aut\Gamma$ has a cyclic subgroup  that is regular or bi-regular on $E$ with $t$ orbits on $V$, and let $\lambda$ be a positive integer. Then, for 
the $\lambda$-extender $\Ga^{(\lambda)}$,
\begin{enumerate}
\item[(a)]  $\Aut \Ga^{(\lambda)}$ has a cyclic regular or bi-regular subgroup, respectively, with $t$ orbits on $V$; 
\item[(b)] if $C=(e_1,\dots, e_{\ell})$ is an Euler cycle of $\Gamma$, and if $\psi\in\Aut \Ga$ induces one of $\varphi, \varphi^2, \tau, \varphi\tau$ in $H(C)$, as in \eqref{e:varphi} and \eqref{e:varphitau}, then 
\[
C^{(\lambda)}=(e_1^1,\dots,e_{\ell}^1,e_1^2,\dots,e_{\ell}^2,\dots,e_1^\lambda,\dots,e_{\ell}^\lambda)
\]
is an Euler cycle of $\Gamma^{(\lambda)}$, and the map $\psi^{(\lambda)}$ as in Table~\ref{t:psilambda} lies in $\Aut\, \Gamma^{(\lambda)}$ and induces the element $\varphi, \varphi^2, \tau, \varphi\tau$ in $H(C^{(\lambda)})$, respectively.
In particular, if $\Ga$ has a symmetrical Euler cycle, 
then also $\Ga^{(\lambda)}$ has a symmetrical Euler cycle.
\end{enumerate}
\end{proposition}

\begin{table}
\begin{center}
\caption{Elements in $H(C^{(\lambda)})$ for an Euler cycle $C$ of $\Gamma$; in all cases $\psi^{(\lambda)}|_V=\psi|_V$ }\label{t:psilambda}
\begin{tabular}{|c|l|}
\hline 
$\psi$ induces in $H(C)$ & $\psi^{(\lambda)}$ on edges, for $i$ mod $\ell$, $j$ mod  $\lambda$   \\  \hline 
$\varphi$ &  $e_i^j\to \left\{\begin{array}{lll}
			e_{i+1}^j &\ &\mbox{if $1\leq i\leq \ell-1$}\\
			e_{i+1}^{j+1} &\ &\mbox{if $i= \ell$}\\
				\end{array}\right. $  \\  
\hline 
$\varphi^2$     &  $e_i^j\to \left\{\begin{array}{lll}
			e_{i+2}^j &\ &\mbox{if $1\leq i\leq \ell-2$}\\
			e_{i+2}^{j+1} &\ &\mbox{if $i\in\{ \ell-1, \ell\}$}\\
				\end{array}\right. $\\  
\hline 
$\tau$ &  $e_i^j\to \left\{\begin{array}{lll}
			e_{\ell+1-i}^{\lambda+1-j} &\ &\mbox{if $1\leq i\leq \ell-1$}\\
			e_{\ell+1-i}^{\lambda-j} &\ &\mbox{if $i= \ell$}\\
				\end{array}\right. $  \\  
\hline 
$\varphi\tau$ &  $e_i^j\to \left\{\begin{array}{lll}
			e_{\ell-i}^{\lambda+1-j} &\ &\mbox{if $1\leq i\leq \ell-2$}\\
			e_{\ell-i}^{\lambda-j} &\ &\mbox{if $i= \ell-1$}\\
			e_{\ell}^{\lambda-j} &\ &\mbox{if $i= \ell$}\\
				\end{array}\right. $  \\  
\hline 
\end{tabular} 
\end{center}
\end{table}

\proof
(a) By assumption, there is an element  $\psi\in\Aut\Ga$ such that $G=\l \psi\r$ is regular or bi-regular on $E$. Thus we may label the edge-set $E=\{e_1,\dots, e_{\ell}\}$, 
such that $e_i=[\b_{i},e_i,\a_i]$ for each $i$, and for $1\leqslant i\leqslant\ell$, reading subscripts modulo $\ell$,
\[
\psi:\ e_{i}\to e_{i+1},\  \mbox{if $\psi$ is regular on $E$, or}\ \psi:\ e_{i}\to e_{i+2},\  \mbox{if $\psi$ is bi-regular on $E$,}
\]
and hence  $\psi:\{\b_{i},\a_i\}\to \{\b_{i+1},\a_{i+1}\}$ or $\{\b_{i},\a_i\}\to \{\b_{i+2},\a_{i+2}\}$, respectively, for each $i$.
%
For each $i$, we  label the $\lambda$ edges of $\Ga^{(\lambda)}=(V,E^{(\lambda)})$ which correspond to the edge  $e_i=[\b_{i},e_i,\a_i]$ by $e_i^j$, for $1\leqslant j\leqslant \lambda$, as in \eqref{e:Elambda}, 
and we define the map $\psi^{(\lambda)}: V\cup E^{(\lambda)}\to V\cup E^{(\lambda)}$ (reading subscripts modulo $\ell$, and  superscripts modulo $\lambda$) by 
$\rho|_V=\psi|_V$ and
\[
\psi^{(\lambda)}: e_i^j\to \left\{\begin{array}{lll}
			e_{i+1}^j &\ &\mbox{if $1\leq i\leq \ell-1$}\\
			e_{i+1}^{j+1} &\ &\mbox{if $i= \ell$}\\
				\end{array}\right.
\]
if $\l\psi\r$ is regular, or 
\[
\psi^{(\lambda)}: e_i^j\to \left\{\begin{array}{lll}
			e_{i+2}^j &\ &\mbox{if $1\leq i\leq \ell-2$}\\
			e_{i+2}^{j+1} &\ &\mbox{if $i\in\{ \ell-1, \ell\}$}\\
				\end{array}\right.
\]
if $\l\psi\r$ is bi-regular. Then $\psi^{(\lambda)}$  is a bijection and preserves incidence in  
$\Ga^{(\lambda)}$, and hence $\psi^{(\lambda)}\in\Aut \Ga^{(\lambda)}$. Further $\l\psi^{(\lambda)}\r$ is regular or bi-regular on $E^{(\lambda)}$, according as $G$ is regular or bi-regular on $E$, respectively, and as $\psi^{(\lambda)}|_V=\psi|_V$, the groups $\l\psi\r, \l\psi^{(\lambda)}\r$ have the same number of vertex-orbits. This proves part (a).

(b) Suppose now that $\Ga$ has an Euler cycle $C=(e_1,\dots, e_{\ell})$ with each $e_i=[\a_{i-1},e_i,\a_i]$, and that $\psi$ induces one of the maps $\varphi, \varphi^2, \tau, \varphi\tau$ in $H(C)$. Let $C^{(\lambda)}$ and $\psi^{(\lambda)}$ be as in the statement and Table~\ref{t:psilambda}. From the definition of $E^{(\lambda)}$ it is clear that $C^{(\lambda)}$ is a cycle, and hence an Euler cycle of $\Gamma^{(\lambda)}$. Clearly $\psi^{(\lambda)}$ is a bijection on $V\cup E^{(\lambda)}$ and preserves incidence and hence lies in $\Aut\,\Gamma^{(\lambda)}$. A somewhat tedious checking shows that the map 
$\psi^{(\lambda)}$ is equal to the element $\varphi, \varphi^2, \tau, \varphi\tau$ in $H(C^{(\lambda)})$, respectively.  In particular if $C$ is a symmetrical Euler cycle of $\Gamma$, so there is an element $\psi\in\Aut \Gamma$ inducing $\varphi^2$ in $H(C)$, then we have just shown that $\psi^{(\lambda)}$ induces $\varphi^2$ in $H(C^{(\lambda)})$. Hence $C^{(\lambda)}$ is a symmetrical Euler cycle in $\Ga^{(\lambda)}$.
  \qed





\section{Examples of graphs}\label{sec:examples}

Here we introduce the families of graphs occurring in the tables for Theorem~\ref{main-thm}. In the light of Proposition~\ref{p:biregext}, to show that a $\lambda$-extender $\Gamma=\Gamma_0^{(\lambda)}$ satisfies the hypotheses of Theorem~\ref{main-thm}, it is sufficient to show that the graph $\Gamma_0$ satisfies them.  Thus we will in general consider graphs which are not $\lambda$-extenders.


\subsection{Circulants}\label{s:circ}

View $\ZZ_n=\{0,1,\dots,n-1\}$ as a cyclic group of order $n$ under addition. The circulants we construct may have multiple edges. 
Let $S$ be a multiset  of elements from $\ZZ_n\setminus\{0\}$ such that, if an element $z\in \ZZ_n$ 
appears exactly $\lambda$ times in $S$, we write $z^{(\lambda)}\in S$.
Assume further that $S$ is \emph{self-inverse}, that is, if $z^{(\lambda)}\in S$ then 
$(-z)^{(\lambda)}\in S$. For such a self-inverse multiset $S$, $\Circ(n,S)$ denotes the (Cayley) graph with vertex set $\ZZ_n$ such that, for each $v\in \ZZ_n$ and each $z^{(\lambda)}\in S$, there are exactly $\lambda$ edges between $v$ and $v+z$. Graphs of this form are called {\it circulants of order $n$}. If each element of a multiset has the same multiplicity $\lambda$, then we sometimes denote the multiset by $S^{(\lambda)}$ for some subset $S\subseteq \ZZ_n$, that is, $S^{(\lambda)}=\{ s^{(\lambda)}\mid s\in S\}$.  First we construct the graphs occurring in Theorem~\ref{main-thm} which are essentially cycles.

\begin{definition}\label{n-cycle}
{\rm
Let $n$ be an integer, $n\geq3$, and let $S=\{1,-1\}$, regarded as  a multiset in $\ZZ_n$. Let $\Gamma=\Circ(n,S)=(V,E)$ with $V=\ZZ_n$ and $E=\{e_i\mid i\in\ZZ_n\}$, where $e_{i}=[i,e_{i},i+1]$, for $i\in\ZZ_n$. Define the map $g:V\cup E\to V\cup E$ by 
\[
g: i\to i+1,\ \mbox{and }\  g: e_{i}\to e_{i+1}\ \mbox{ for $i\in\ZZ_n$}. 
\]
Define the edge sequence $C=(e_{0},\dots,e_{n})$ of length $n$.

}
\end{definition}

\begin{lemma}\label{l:n-cycle}
With the notation of Definition~$\ref{n-cycle}$, the graph $\Gamma\cong \C_n$, $\Gamma$ is connected,  $g\in\Aut\Gamma$, and 
\begin{enumerate}
\item[(a)] the group $\la g\ra$ is transitive on $V$ and regular on $E$ and line $1$ of Table~\ref{t:mainreg} holds. Moreover, $C$ is a symmetrical Euler cycle for $\Gamma$, $g$ 
induces the map $\varphi$ of \eqref{e:varphi}, and $H(C)=D(C)\cong \D_{2n}$.

\item[(b)] For the $2$-extender $\C_n^{(2)}$, the group $\la g\ra$ is transitive on $V$ and bi-regular on $E^{(2)}$ and line $1$ of Table~\ref{t:mainbireg1} holds. 

\item[(c)] If $n$ is even, $n\geq4$,  then $\l g^2\r\cong \ZZ_{n/2}$ is bi-transitive on $V$, and bi-regular on $E$ with orbits $E_0$ and $E_1$ such that $[E_0]\cong [E_1]\cong \frac{n}{2}\K_2$, and line $1$ of Table~\ref{t:mainbireg2} holds. 
\end{enumerate}
\end{lemma}

\proof
Clearly $\Gamma\cong \C_n$, $\Gamma$ is connected, and  $g\in\Aut\Gamma$, and most of part (a) follows immediately from the definitions of $g$ and a symmetrical Euler cycle for $\Gamma$,  with $g$ inducing the map $\varphi$ in $H(C)$.
Further $\Aut\Gamma$ contains an automorphism $\tau$ as in \eqref{e:varphi}, and $\l g, \tau\r\cong \D_{2n}$, so $H(C)=D(C)=\D_{2n}$. For part (b), we note that $\l g\r$ has two orbits in $E^{(2)}$, namely, for $j\in\{1,2\}$, $E_j:=\{e_i^j\mid i\in\ZZ_n\}$, where $e_{i}^j=[i,e_{i}^j,i+1]$ for $i\in\ZZ_n$. 
Finally for part (c), for $n$ even, $\l g^2\r\cong \ZZ_{n/2}$ is bi-transitive on $V$, and bi-regular on $E$ with edge-orbits $E_0:=\{e_{2i}\mid 0\leq i<n/2\}$ and $E_1:=\{e_{2i+1}\mid 0\leq i<n/2\}$, and the induced subgraphs $[E_0]\cong [E_1]\cong \frac{n}{2}\K_2$.
\qed

Next we consider the remaining circulants occurring in Table~\ref{t:mainbireg1} for Theorem~\ref{main-thm}. For $a\in\ZZ_n$, by $|a|$ we mean the additive order of $a$, that is, the least positive integer $m$ such that $ma\equiv 1\pmod{n}$.

\begin{definition}\label{d:circ1}
{\rm
Let $n$ be a positive integer with $n\geq3$,  and let $a,b\in \ZZ_n\setminus\{0\}$ such that  $\gcd(n,a,b)=1$.  Let $\Ga=\Ga(n,a,b)=\Circ(n,S)$ where 
\[
S= \left\{\begin{array}{lll}
			\{a,-a,b,-b\} &\ &\mbox{if $|a|\geq3$, $|b|\geq3$, $a\ne\pm b$}\\
			\{a,-a,b^{(2)}\} &\ &\mbox{if $|a|\geq3$, $2b=n$}\\
				\end{array}\right.
\]
so $\Ga$ has $2n$ edges, namely $e_{i,a}=[i,e_{i,a}, i+a]$, and $e_{i,b}=[i,e_{i,b}, i+b]$, for $i\in\ZZ_n$, and we note that, if $2b=n$ then $e_{i,b}$ and $e_{i+b,b}$ are both incident with the same pair of vertices $\{i,i+b\}$. Let     
$
E_a=\{e_{i,a}\mid i\in\ZZ_n\},\ E_b=\{e_{i,b}\mid i\in\ZZ_n\}. 
$
Thus $\Gamma$ is not always simple, but in all cases the natural generator $g$ of $\ZZ_n$ acts as follows:
\[
g: i\to i+1,\quad e_{i,a}\to e_{i+1,a},\quad e_{i,b}\to e_{i+1, b},\ \mbox{for $i\in\ZZ_n$.}
\] 
}
\end{definition}

\begin{lemma}\label{l:circ1}
With the notation of Definition~$\ref{d:circ1}$, the graph $\Gamma=\Gamma(n,a,b)$ is connected, and the group  
$\l g\r=\ZZ_n$ induces a cyclic subgroup of $\Aut\Ga$ which is regular  on vertices and bi-regular on edges, with edge-orbits $E_a$ and $E_b$. Moreover,  if $2b\ne n$ then $\Gamma$ is  simple of valency $4$ as in line $4$ of Table~\ref{t:mainbireg1}.
\end{lemma}

\proof
The vertices reached by sequences of edges beginning at $0$ are precisely those of the form $ia+jb$ for some integers $i, j$, and these vertices are the multiples of $\gcd(a,b)$, modulo $n$. Since $\gcd(n,a,b)=1$ it follows that all vertices occur so $\Gamma$ is connected and  $\ZZ_n=\l a,b\r$. By definition of circulants, $\ZZ_n$, acting naturally by `right multiplication' induces a cyclic subgroup of $\Aut\Ga$ which is regular  on vertices, and it has $E_a, E_b$ as its edge-orbits. If $2b\ne n$ then $\Gamma$ has valency $|S|=4$ and is a graph in  line $4$ of Table~\ref{t:mainbireg1}.
\qed

Note that if $|b|=2$ in Definition~\ref{d:circ1} then $n$ is even and $b=n/2$.   
The following corollary describes two special cases of the construction in Definition~\ref{d:circ1} where $|b|=2$. We note that in case (b) below, the graph $\Gamma(2r,2,r)$ is isomorphic to the cartesian product $\C_r\Box\K_2^{(2)}$.

\begin{corollary}\label{c:circ1}
Let $n=2r\geq 4$. Then the following hold, with the notation of 
Definition~$\ref{d:circ1}$,
\begin{enumerate}
\item[(a)] $\Gamma=\Gamma(2r,1,r) = \Circ(2r,\{1,-1\}) + \Circ(2r,\{r^{(2)}\}) 
= \C_{2r} + r\K_2^{(2)}$; and
\item[(b)] if $r$ is odd, then $\Gamma=\Gamma(2r,2,r) = \Circ(2r,\{2,-2\}) + \Circ(2r,\{r^{(2)}\}) = 2\C_{r} + r\K_2^{(2)}$.
\end{enumerate}
In either case, $\Gamma$ is connected, and $\Aut\Gamma$ has a cyclic subgroup that is transitive on vertices and bi-regular on edges, as in line $2$ or $3$ of Table~\ref{t:mainbireg1}.
\end{corollary}

\proof
This follows immediately from Proposition~\ref{l:circ1}, noting that $\gcd(2r,1,r)=1$, and when $r$ is odd also $\gcd(2r,2,r)=1$. 
\qed

\subsection{Modified Kronecker product graphs}
Our next construction is a modification of the Kronecker product construction for graphs, so we use similar notation. It produces the graphs in  line $3$ of Table~\ref{t:mainbireg2}.

\begin{definition}\label{d:cyclest}
{\rm
Let $r,s,t$ be positive integers such that $r\geq2$,  $st\geq2$ and $\gcd(s,t)=1$,  and let $S=\ZZ_s, T=\ZZ_t$. Then $\Ga=\C_{2r}[s\K_1, t\K_1]=(V,E)$ is the graph defined as follows. The vertex set $V= V_S\cup V_T$, where
\[
V_S = \{ (2k,i)\mid 0\leq k\leq r-1,\ i\in S\}, \mbox{and}\ 
V_T=\ \{ (2k+1,j)\mid 0\leq k\leq r-1,\ j\in T\},
\]
and edge set $E= E_S\cup E_T\ \mbox{where}\ $
\begin{align*}
E_S &= \{ [(2k,i), e_{i,j}^{2k}, (2k+1,j)] \mid 0\leq k\leq r-1, i\in S, j\in T\}, \mbox{and}\\
E_T &= \{ [(2k+1,j), e_{j,i}^{2k+1}, (2k+2,i)] \mid 0\leq k\leq r-1, i\in S, j\in T\},
\end{align*}
and we read the entry $2k+2$ modulo $2r$. 
Define the maps  $g, y:V\cup E \to V\cup E$ as follows, for $0\leq k\leq r-1, i\in S, j\in T$:
\[
\begin{array}{lclcclcl}
g: 	&(2k,i)&\to &(2k+2,i), &(2k+1,j)	&\to &(2k+3,j)&\ \mbox{if $k\leq r-2$,}\\ 
	&(2r-2,i)&\to &(0,i+1),  &(2r-1,j)	&\to &(1,j+1) & \\
	&e_{i,j}^{2k}&\to &e_{i,j}^{2k+2},\ &e_{j,i}^{2k+1}&\to &e_{j,i}^{2k+3} &\ \mbox{if $k\leq r-3$,}\\ 
	&e_{i,j}^{2r-4}&\to &e_{i,j}^{2r-2},\ &e_{j,i}^{2r-3}&\to &e_{j,i+1}^{2r-1}\  &\\ 
	&e_{i,j}^{2r-2}&\to &e_{i+1,j+1}^{0},\ &e_{j,i}^{2r-1}&\to &e_{j+1,i}^{1}\  &\\
	\end{array}
\]
\[
\begin{array}{lclcclcl}
y: 	&(0,i)&\to &(0, -i), &(1,j)	&\to &(2r-1, -j-1)&\ \\ 
 	&(2k,i)&\to &(2(r-k) , -i-1), &(2k+1,j)	&\to &(2(r-k-1)+1, -j-1)&\ \\ 
	&e_{i,j}^{0}&\to &e_{-j-1,-i}^{2r-1},\ &e_{j,i}^{2r-1}&\to &e_{-i,-j-1}^{0} &\\ 
	&e_{i,j}^{2k}&\to &e_{-j-1,-i-1}^{2(r-k-1)+1},\ &e_{j,i}^{2k-1}&\to &e_{-i-1,-j-1}^{2(r-k)} &\\ 
\end{array}
\]
where $i\in S, j\in T$, and in lines 2 and 4 in the defintion of $y$ we have $1\leq k\leq r-1$.
}
\end{definition}

We could have defined the graph $\C_{2r}[s\K_1, t\K_1]$ with $s=t=1$, but then it is just the cycle $\C_{2r}$ so we avoid this degenerate case.
It is not difficult to see that the automorphism group of $\C_{2r}[s\K_1, t\K_1]$ is isomorphic to $(\Sym(s)\times\Sym(t))\wr\D_{2r}$.
We show that the map $g$ generates a vertex-bi-transitive, edge-bi-regular group of automorphisms, and moreover  that $\l g,y\r$ is an edge-regular dihedral group.

\begin{center}
\begin{figure}
 \begin{tikzpicture}
      \foreach \k in {0,2,4,6} { 
        \foreach \i in {0,1} {
        \coordinate (v\k\i) at ($(120-\k*60: 3.5cm-1cm*\i)$);
        };
      };
      \foreach \k in {1,3,5} { 
        \foreach \i in {0,1,2} {
        \coordinate (v\k\i) at ($(120-\k*60: 3.5cm-0.5cm*\i)$);
        }; 
      };

     \draw[rounded corners=6pt, thick, color=blue!50, rotate =60] (2,-.5) rectangle (4,.5);
     \draw[rounded corners=6pt, thick, color=blue!50, rotate =180] (2,-.5) rectangle (4,.5);
     \draw[rounded corners=6pt, thick, color=blue!50, rotate =300] (2,-.5) rectangle (4,.5);
     \draw[rounded corners=6pt, thick, color=red!50] (2,-.5) rectangle (4,.5);
     \draw[rounded corners=6pt, thick, color=red!50, rotate =120] (2,-.5) rectangle (4,.5);
     \draw[rounded corners=6pt, thick, color=red!50, rotate =240] (2,-.5) rectangle (4,.5);

      \foreach \k/\l in {0/1,2/3,4/5} {
      \foreach \i  in {0,1} {
        \foreach \j in {0,1,2} {
         \draw[color=red!50] (v\k\i)--(v\l\j); 
        };
      };
      };

      \foreach \k/\l in {1/2,3/4,5/0} {
      \foreach \i  in {0,1} {
        \foreach \j in {0,1,2} {
         \draw[color=blue!50] (v\k\j)--(v\l\i); 
        };
      };
      };

      \draw[line width= 1.5pt, color =red] 
          (v00)-- node[above] {$e_{0,0}^{0}$}(v10) 
          (v20)-- node[right] {$e_{0,0}^{2}$}(v30) 
          (v40)-- node[left] {$e_{0,0}^{4}$}(v50) 
          (v01)-- node[below] {$e_{1,1}^{0}$}(v11) 
          (v41)-- node[right] {$e_{1,2}^{4}$}(v52); 
      \draw[line width= 1.5pt, color =blue] 
          (v10)-- node[right] {$e_{0,0}^{1}$}(v20) 
          (v30)--node[below] {$e_{0,0}^{3}$}(v40) 
          (v50)--node[left, near start] {$e_{0,1}^{5}$}(v01) 
          (v52)--node[left, near end] {$e_{2,0}^{5}$}(v00) 
          (v11)--node[left] {$e_{1,1}^{1}$}(v21); 

      \foreach \k in {0,2,4} {\foreach \i in {0,1} { \fill (v\k\i) circle (2pt);}; };
      \foreach \k in {1,3,5} {\foreach \i in {0,1,2} { \fill (v\k\i) circle (2pt);}; };
      \node[anchor =south east] at (v00) {\tiny $(0,0)$};
      \node[anchor =north west] at (v01) {\tiny $(0,1)$};
      \node[anchor =south west] at (v10) {\tiny $(1,0)$};
      \node[anchor = east] at (v11) {\tiny $(1,1)$};
      \node[anchor =north east] at (v12) {\tiny $(1,2)$};
      \node[anchor =west] at (v20) {\tiny $(2,0)$};
      \node[anchor =east] at (v21) {\tiny $(2,1)$};
      \node[anchor =north west] at (v30) {\tiny $(3,0)$};
      \node[anchor =south east] at (v32) {\tiny $(3,2)$};
      \node[anchor =north east] at (v40) {\tiny $(4,0)$};
      \node[anchor =south west] at (v41) {\tiny $(4,1)$};
      \node[anchor =east] at (v50) {\tiny $(5,0)$};
      \node[anchor =west] at (v52) {\tiny $(5,2)$};
 \end{tikzpicture}
\caption{
Graph $\Ga=\C_{2r}[s\K_{1},t\K_1]$ from Definition~\ref{d:cyclest}, with $r=3, s=2, t=3$
}\label{fig-crst}
\end{figure}
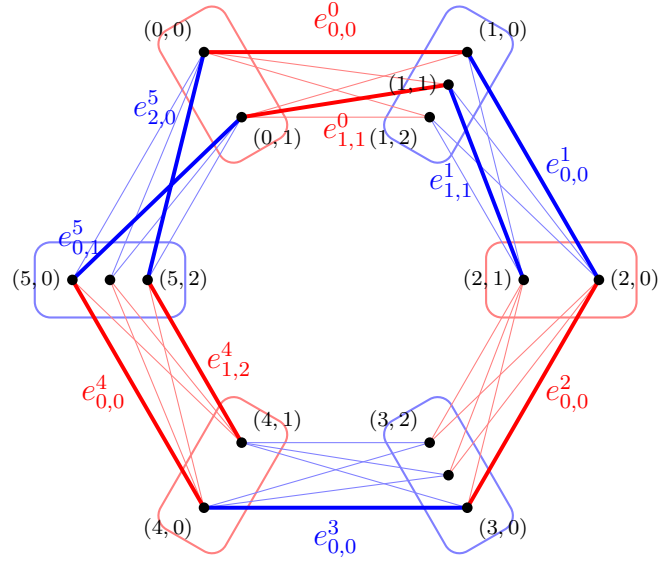
\end{center}

\begin{lemma}\label{l:cyclest}
With the notation of Definition~$\ref{d:cyclest}$, $\Ga=\C_{2r}[s\K_1, t\K_1]$ is connected and is the graph in line $3$ of Table~\ref{t:mainbireg2}, the maps $g, y$ are automorphisms of 
$\Gamma$, and  
\begin{enumerate}
\item[(a)]
$\l g\r $ is cyclic of order $rst$, bi-transitive on $V$ with orbits $V_S, V_T$, and bi-regular on $E$ with orbits $E_S, E_T$; also $\l g,y\r \cong \D_{2rst}$ and is regular on $E$;
\item[(b)] further, for each positive integer $\lambda$, the graph $\Gamma^{(\lambda)}$,  admits a cyclic subgroup $\ZZ_{rst\lambda}$, and a dihedral subgroup $\D_{2rst\lambda}$ of automorphisms that are both bi-transitive on vertices, and are bi-regular or regular on edges, respectively.
\end{enumerate}
\end{lemma}

\proof
It is straightforward to check that $\Gamma$ is connected, and that $\Gamma$ is the graph in line $3$ of Table~\ref{t:mainbireg2}. Also $g, y$ are bijections and it is straightforward, if tedious, to check that each of $g, y$ preserves incidence. Hence $g,y\in\Aut\Gamma$. First we consider repeated applications of $g$ on the edge $e=e_{0,0}^0=[(0,0),e_{0,0}^0,(1,0)]\in E_S$. Suppose that $e^{g^m}=e$. Now $e^{g^m}=e^{2m}_{i,j}$ for some $i\in S,j\in T$, reading the superscript $2m$ modulo $2r$, and this means, since $e^{g^m}=e$, that $m=r\ell$ for some integer 
$\ell$. Now $e^{g^r}=  [(0,1), e_{1,1}^{0}, (1,1)]$, and repeating this we find that $e^{g^{r\ell}} = [(0,\ell), e_{\ell,\ell}^{0}, (1,\ell)]$, where in the first entry we read $\ell$ modulo $s$, and in the third entry we read $\ell$ modulo $t$. Thus, again since $e^{g^{r\ell}}=e^{g^m}=e$, we conclude that $\ell$ is divisible by both $s$ and $t$, and hence by $st$ since $\gcd(s,t)=1$.  Thus $m$ is a multiple of $rst$, and since $|E_S|=rst$ and $g$ leaves $E_S$ invariant it follows that the $\l g \r$-orbit containing $e$ is $E_S$. A similar proof shows that the  $\l g \r$-orbit containing $e^1_{0,0}$ is $E_T$, and it follows that $\l g\r$ is bi-regular on $E$. Also the vertex-orbits of $\l g\r$ are clearly $V_S$ and $V_T$ and we conclude that $\l g\r $ has order $rst$, and is bi-transitive on $V$. This proves the first assertion.

Now we consider $y$. It follows from the previous paragraph that $y^2, (yg)^2\in \Aut\Gamma$. Again it is tedious, but not too difficult, to verify that both $y^2$ and $(yg)^2$ act trivially on the vertex set $V$, and hence each of them leaves invariant the edge subsets $E_{S,k}=\{ e^{2k}_{i,j}\mid i\in S, j\in T\}$ and $E_{T,k}=\{ e^{2k+1}_{j,i}\mid i\in S, j\in T\}$, for each $k=0,\dots, r-1$. For a fixed $k$, since each vertex pair $\{(2k,i), (2k+1,j)\}$ (for $i\in S, j\in T$) is incident with a unique edge $e^{2k}_{i,j}\in E_{S,k}$, it follows that $y^2$ and $(yg)^2$ act trivially on $E_{S,k}$. Similarly $y^2$ and $(yg)^2$ act trivially on $E_{T,k}$, and since this holds for each $k$, it follows that $y^2=(yg)^2=1$. Thus $ygy=g^{-1}$ and so $\l g,y\r\cong\D_{2rst}$. Moreover, since $y$ interchanges $E_S$ and $E_T$ it follows that $\l g,y\r$ is regular on the edge set $E$. This completes the proof of part (a).
Part (b) now follows immediately from Proposition~\ref{p:biregext}(b). 
\qed

\subsection{Complete bipartite graphs}\label{s:kst}

Next we study complete bipartite graphs, their $\lambda$-extenders, and constructions which combine complete bipartite graphs with circulants. 
First we establish the properties of these graphs occurring in  Tables~\ref{t:mainreg} 
and~\ref{t:mainbireg2}.

\begin{definition}\label{d:kst}
{\rm
Let $s, t$ be positive integers such that  $st\geq2$ and $\gcd(s,t)=1$,  and let $S=\ZZ_s, T=\ZZ_t$. Then $\Gamma=\K_{s,t} = (V, E)$ with $V=S\cup T$ and 
\[
E= \{ [i, e_{i,j}, j] \mid i\in S, j\in T\}.
\]
Define the map $g:V\cup E \to V\cup E$ as follows, for $i\in S, j\in T$:
\[
g: i\to i+1 \pmod{s},\quad  j\to j+1 \pmod{t},\quad e_{i,j}\to e_{i+1,j+1}.
\] 
}
\end{definition}

\begin{lemma}\label{l:kst}
With the notation of Definition~$\ref{d:kst}$ so in particular $\gcd(s,t)=1$, 
the map $g$ is an automorphism of 
$\Gamma=\K(s,t)$, and  $\l g\r \cong \ZZ_{st}$ is cyclic, bi-transitive on $V$ with orbits $S, T$, and regular on $E$, and $\Ga$ is as in line $2$ of Table~\ref{t:mainreg}. Further, the graph $\K(s,t)^{(2)}$ admits a cyclic subgroup of automorphisms that is bi-transitive on vertices, and  bi-regular on edges, as in line $2$ of Table~\ref{t:mainbireg2}. 
\end{lemma}

\proof
By definition both  $g|_V$ and  $g|_E$ are bijections, and $g$ preserves incidence, so $g\in\Aut\Gamma$. Also, since $\gcd(s,t)=1$, $|g|=st$, and it follows that the cyclic group $\l g\r$ is bi-transitive on $V$ with orbits $S, T$, and regular on $E$. The fact that 
$\Aut  \K(s,t)^{(2)}$ has a cyclic subgroup bi-transitive on vertices, and bi-regular  on edges follows on considering $\K(s,t)^{(2)}$ as an edge-disjoint union of two copies of $\K(s,t)$ on the same vertex set.
\qed

The next construction starts with a complete bipartite graph and adds a second set of edges incident with the vertices of one bipart which form a perfect matching or an edge-disjoint union of cycles (with multiple edges). These are the graphs in lines $4$ and $6$ of Table~\ref{t:mainbireg2}.

\begin{definition}\label{d:ck}
{\rm
Let $r,n,t$ be positive integers such that $rt\geq1$, $n\geq2$, and $\gcd(nr,t)=1$,  
and let $V_1=\ZZ_{rn}, V_2=\ZZ_t$. Define a graph $\Ga=\C\K(r,n,t)=(V,E)$ with 
vertex-set $V= V_1\cup V_2$, and edge-set $E=E_0\cup E_1$ where
\begin{align*}
E_0 &= \{ [k, e_{k,j}, j] \mid k\in V_1, j\in V_2\},\quad \mbox{and}\\
E_1 &= \{ [k, e_{k}^{j}, k+r]\mid k\in V_1, j\in V_2\}.
\end{align*}
Define the map  $g:V\cup E \to V\cup E$ as follows, for $k\in V_1, j\in V_2$:
\[
\begin{array}{lclcclcl}
g: 	&k&\to &k+1, &j	&\to &j+1&\ \\ 
	&e_{k,j}&\to &e_{k+1,j+1},\ &e_{k}^{j}&\to &e_{k+1}^{j+1} &\ \\ 
\end{array}
\]
}
\end{definition}

\begin{lemma}\label{l:ck}
With the notation of Definition~$\ref{d:ck}$ so in particular $\gcd(rn,t)=1$, 
the graph $\Gamma=\C\K(r,n,t)$ is connected and the map $g$ is an automorphism of 
$\Gamma$, and  $\l g\r \cong \ZZ_{rnt}$ is cyclic, bi-transitive on $V$ with orbits $V_1, V_2$, and bi-regular on $E$ with orbits $E_0, E_1$. Further $\Gamma= [E_0]+[E_1]$, where $[E_0] = \K_{rn,t}$ and, if $n=2$ then $[E_1]= r\K_2^{(2t)}$  as in line $4$ of Table~\ref{t:mainbireg2}, while if $n\geq3$ then $[E_1]= r\C_n^{(t)}$  as in line $6$ of Table~\ref{t:mainbireg2}.
\end{lemma}

\proof
It is clear that $\Gamma= [E_0]+[E_1]$ is connected and $[E_0] = \K_{rn,t}$. Further, if $n=2$ then for each  $k\in V_1, j\in V_2$, both $e_{k}^j$ and $e_{k+r}^j$ are incident with $k$ and $k+r$, and hence $[E_1]= r\K_2^{(2t)}$. On the other hand if $n\geq3$ then $[E_1]= r\C_n^{(t)}$, so the last assertions hold.
  
By definition, both  $g|_V$ and  $g|_E$ are bijections, and $g$ preserves incidence, so $g\in\Aut\Gamma$. Consider the action of $g$ on $E$. First let $e=e_0^0\in E_1$, and suppose that $m$ is the least positive integer such that $e^{g^m}=e$. Now, by definition of $g$, $e^{g^m}=[m,e_m^m,m]$, where the first entry $m\in\ZZ_{rn}$ and the last entry $m\in\ZZ_t$. Thus, since $e^{g^m}=e$ the integer $m$ is divisible by $rn$ and $t$, and hence by $rnt$ since $\gcd(rn,t)=1$. On the other hand $g^{rnt}$ fixes $e$ and hence $m=rnt$ and as $|E_1|=rnt$ we conclude that $\l g\r$ acts regularly on $E_1$. A similar argument shows that $\l g\r$ acts regularly on $E_0$, and so $|g|=rnt$ and $\l g\r$ is bi-transitive on $V$ and bi-regular on $E$.
\qed

In a more general but similar construction, we start with an edge-disjoint union of complete bipartite graphs $u\K_{sr,t}$, and add a second set of edges incident with the vertices of the bipart of size $usr$, which form either a perfect matching of this bipart, or an edge-disjoint union of cycles. In this case no edge of the second edge-set is incident with two vertices of the same component $\K_{sr,t}$. These are the graphs in lines $5$ and $7$ of Table~\ref{t:mainbireg2}, and are depicted in Figure~\ref{fig-ck}.

\begin{definition}\label{d:ck2}
{\rm
Let $r,s,t,u$ be positive integers such that  $u\geq2$, and $\gcd(r,u)=\gcd(sr,t)=1$,  
and let $V_1=\ZZ_{sru}, V_2=\ZZ_{ut}$. Define a graph $\Ga=\C\K_{(2)}(r,s,t,u)=(V,E)$ with 
vertex-set $V= V_1\cup V_2$, and edge-set $E=E_0\cup E_1$ where
\begin{align*}
E_0 &= \{ [i+uk, e_{i, k ,j}, i+uj] \mid 0\leq i\leq u-1, 0\leq k\leq sr-1, 0\leq j\leq t-1\},\quad \mbox{and}\\
E_1 &= \{ [\ell + rk, e_{\ell,k}^{j}, \ell + r(k+1)]\mid 0\leq \ell\leq r-1,  0\leq k\leq su-1, 0\leq j\leq t-1\} .
\end{align*}
Note that for edges in $E_0$,  we read the first entry as an element of $\ZZ_{sru}$ and the third entry as an element of $\ZZ_{ut}$, while for edges in $E_1$, the first and third entries are elements of $\ZZ_{sru}$.
Define the map  $g:V\cup E \to V\cup E$ as follows, for $k\in V_1, j\in V_2$:
\[
\begin{array}{lcllclcl}
g: 	&k&\to &k+1, \ \mbox{and } j \to j+1, &	& &&\ \\ 
	&e_{i,k,j}&\to &e_{i+1,k,j+1},\ &  & &\ \\ 
	&e_{\ell,k}^{j}&\to &e_{\ell+1,k}^{j+1},\ \mbox{if $0\leq \ell\leq r-2$,\ and }&e_{r-1,k}^{j}&\to &e_{0,k+1}^{j+1} &\ \\ 
\end{array}
\]
where in the edge-image $e_{i+1,k,j+1}$ we read $i+1$ modulo $u$ and $j+1$ modulo $t$, and for the edge-image $e_{\ell+1,k}^{j+1}$ we read $\ell+1$ modulo $r$ and $j+1$ modulo $t$. 
}
\end{definition}

\begin{lemma}\label{l:ck2}
With the notation of Definition~$\ref{d:ck2}$,  the graph $\Gamma=\C\K_{(2)}(r,s,t,u)$ is connected and the map $g$ is an automorphism of $\Gamma$, and  $\l g\r \cong \ZZ_{srut}$ is cyclic, bi-transitive on $V$ with orbits $V_1, V_2$, and bi-regular on $E$ with orbits $E_0, E_1$. Further $\Gamma= [E_0]+[E_1]$, where $[E_0] = u\K_{sr,t}$ and, 
if $(s,u)=(1,2)$ then $[E_1]= r\K_2^{(2t)}$  as in line $5$ of Table~\ref{t:mainbireg2}, and otherwise $su\geq3$ and $[E_1]= r\C_{su}^{(t)}$  as in line $7$ of Table~\ref{t:mainbireg2}.
\end{lemma}

\proof
By definition $\Gamma= [E_0]+[E_1]$ and $[E_0] = u\K_{sr,t}$. Further, if $(s,u)=(1,2)$, then for each  $\ell, j$, both $e_{\ell,0}^j$ and $e_{\ell,1}^j$ are incident with $\ell$ and $\ell+r$ in $V_1$ and hence $[E_1]= r\K_2^{(2t)}$. On the other hand if $(s,u)\ne(1,2)$, then $su\geq3$ and $[E_1]= r\C_{su}^{(t)}$, so the last assertions hold.
  
By definition, both  $g|_V$ and  $g|_E$ are bijections, and $g$ preserves incidence, so $g\in\Aut\Gamma$. Consider the action of $g$ on $E$. First let $e=e_{0,0,0}\in E_0$, and suppose that $m$ is the least positive integer such that $e^{g^m}=e$. By definition of $g$, $e^{g^m}=[m,e_{m,0,m},m]$, where the first entry $m\in\ZZ_{sru}$ and the last entry $m\in\ZZ_{ut}$. Thus, since $e^{g^m}=e$ the integer $m$ is divisible by $sru$ and $ut$, and hence by $\lcm\{sru,ut\}=srut$ since $\gcd(sr,t)=1$. On the other hand $g^{srut}$ fixes $e$ and hence $m=srut$ and as $|E_0|=srut$ we conclude that $\l g\r$ acts regularly on $E_0$. Now let $e=[0,e_{0,0}^0,0]\in E_1$, and suppose that $e^{g^m}=e$ with $m\geq1$ minimal. By definition of $g$, $e^{g^m}= [m,e_{m,k}^m,m]$, for some $k$, where the first entry $m\in\ZZ_{sru}$ and the last entry $m\in\ZZ_{ut}$, and as before we deduce that $m=srut=|E_1|$ and that  $\l g\r$ acts regularly on $E_0$. 
Thus $|g|=srut$ and $\l g\r$ is bi-transitive on $V$ and bi-regular on $E$.
\qed


\begin{remark}\label{r:figck}
We give some details about Figure~\ref{fig-ck},  for the graph $\Gamma=\C\K_{(2)}(r,s,t,u)= r\C_{su}^{(t)}+u\K_{sr,t}$ with $su\geq3$, from Definition~\ref{d:ck2}. 
The subset $V_1$ of vertices admits two $\l g\r$-invariant partitions: first $\{B_1,\dots,B_{r}\}$ where the $B_i$ are the components of the edge-induced subgraph $[E_1]= r\C_{su}^{(t)}$, and second 
 $\{C_1,\dots, C_{u}\}$ where the $C_i$ are the intersections with $V_1$ of the components of $[E_0]= u\K_{sr,t}$. The vertex-subset $V_2$ admits the  $\l g\r$-invariant partition
 $\{D_1,\dots, D_{u}\}$,  where the $D_i$ are the intersections with $V_2$ of the components of $[E_0]= u\K_{sr,t}$.
\end{remark}

\begin{center}
\begin{figure}[ht]

\begin{tikzpicture}
  \draw[rounded corners=6pt, thick, color=blue!50] (0,-.1) rectangle (4,-.9);
  \draw[rounded corners=6pt, thick, color=blue!50] (0,-1.1) rectangle (4,-1.9);
  \draw[rounded corners=6pt, thick, color=blue!50] (0,-3.1) rectangle (4,-3.9);

  \draw[rounded corners=6pt, thick, color=blue!50] (6,-.1) rectangle (8,-.9);
  \draw[rounded corners=6pt, thick, color=blue!50] (6,-1.1) rectangle (8,-1.9);
  \draw[rounded corners=6pt, thick, color=blue!50] (6,-3.1) rectangle (8,-3.9);

  \draw[rounded corners=6pt, thick, color=red!50] (.1,.1) rectangle (.9,-4.1);
  \draw[rounded corners=6pt, thick, color=red!50] (2.1,.1) rectangle (2.9,-4.1);
  \draw[rounded corners=6pt, thick, color=red!50] (3.1,.1) rectangle (3.9,-4.1);

  \draw[decorate,decoration={brace,amplitude=6},thick] (9,-0.3)--(9,-3.7);
  \node at (10,-2) {$u\K_{sr,t} $};
  \draw[decorate,decoration={brace, mirror,amplitude=6},thick] (.5,-5)--(3.5,-5);
  \node at (2,-5.5) {$r\C_{su}^{(t)} $};

  \draw[color=red] (3.3,-.2).. controls (3.1,-.5) and (3.1,-.9) ..(3.3,-1.2);
  \draw[color=red] (3.3,-.2).. controls (3.2,-.5) and (3.2,-.9) ..(3.3,-1.2);
  \draw[color=red] (3.3,-1.2).. controls (3.1,-1.5) and (3.1,-1.9) ..(3.3,-2.2);
  \draw[color=red] (3.3,-1.2).. controls (3.2,-1.5) and (3.2,-1.9) ..(3.3,-2.2);
  \draw[color=red] (3.3,-2.8).. controls (3.1,-3.1) and (3.1,-3.5) ..(3.3,-3.8);
  \draw[color=red] (3.3,-2.8).. controls (3.2,-3.1) and (3.2,-3.5) ..(3.3,-3.8);
  \draw[color=red] (3.3,-3.8).. controls (3.6,-2.1) and (3.6,-2.1) ..(3.5,-.8);
  \draw[color=red] (3.3,-3.8).. controls (3.7,-2.1) and (3.7,-2.1) ..(3.5,-.8);
  \node at (3.3,-2.4) {$\vdots$};
  \node at (2.3,-2.4) {$\vdots$};
  \node at (0.3,-2.4) {$\vdots$};
  \node[color=red!75] at (3.5,-4.5) {\tiny ${B}_1$};
  \node[color=red!75] at (2.5,-4.5) {\tiny ${B}_2$};
  \node[color=red!75] at (1.5,-4.5) { $\cdots$};
  \node[color=red!75] at (.5,-4.5) {\tiny ${B}_{r}$};

  \node[color=blue!75] at (-.5,-.5) {\tiny ${C}_1$};
  \node[color=blue!75] at (-.5,-1.5) {\tiny ${C}_2$};
  \node[color=blue!75] at (-.5,-2.5) { $\vdots$};
  \node[color=blue!75] at (-.5,-3.5) {\tiny ${C}_{u}$};

  \node[color=blue!75] at (8.5,-.5) {\tiny ${D}_1$};
  \node[color=blue!75] at (8.5,-1.5) {\tiny ${D}_2$};
  \node[color=blue!75] at (8.5,-2.5) { $\vdots$};
  \node[color=blue!75] at (8.5,-3.5) {\tiny ${D}_{u}$};

  \coordinate (v111) at (3.3,-0.2);
  \coordinate (v112) at (3.5,-0.8);
  \coordinate (w11) at (6.3, -.2);
  \coordinate (w12) at (6.3, -.5);
  \coordinate (w13) at (6.3, -.8);

  \foreach \y in {1,2,3} {
  \foreach \x in {1,2} {
  \draw [color=blue] (v11\x) .. controls ($0.5*(v11\x) +0.5*(w1\y) +(0,0.1) $) ..(w1\y);
  \draw [color=blue] (v11\x) .. controls ($0.5*(v11\x) +0.5*(w1\y)-(0,0.1) $) ..(w1\y);
  };
  \draw [color=blue] ($(v111)+(-3,0) $) .. controls ($0.5*(v111) +0.5*(w1\y)-(1.5,-1) $) ..(w1\y);
  \draw [color=blue] ($(v111)+(-3,0) $) .. controls ($0.5*(v111) +0.5*(w1\y)-(1.5,-1.3) $) ..(w1\y);
  };

  \fill (3.3, -.2) circle (1.5pt);
  \fill (3.5, -.8) circle (1.5pt);
  \fill (3.3, -1.2) circle (1.5pt);
  \fill (3.3, -2.2) circle (1.5pt);
  \fill (3.3, -2.8) circle (1.5pt);
  \fill (3.3, -3.8) circle (1.5pt);

  \fill (2.3, -.2) circle (1.5pt);
  \fill (2.5, -.8) circle (1.5pt);
  \fill (2.3, -1.2) circle (1.5pt);
  \fill (2.3, -2.2) circle (1.5pt);
  \fill (2.3, -2.8) circle (1.5pt);
  \fill (2.3, -3.8) circle (1.5pt);

  \fill (0.3, -.2) circle (1.5pt);
  \fill (0.5, -.8) circle (1.5pt);
  \fill (0.3, -1.2) circle (1.5pt);
  \fill (0.3, -2.2) circle (1.5pt);
  \fill (0.3, -2.8) circle (1.5pt);
  \fill (0.3, -3.8) circle (1.5pt);

  \fill (6.3, -.2) circle (1.5pt);
  \fill (6.3, -.5) circle (1.5pt);
  \fill (6.3, -.8) circle (1.5pt);
  \node at (7.2,-.5) {$\dots$};

  \fill (6.3, -1.2) circle (1.5pt);1
  \fill (6.3, -1.5) circle (1.5pt);
  \fill (6.3, -1.8) circle (1.5pt);
  \node at (7.2,-1.5) {$\dots$};

  \fill (6.3, -3.2) circle (1.5pt);
  \fill (6.3, -3.5) circle (1.5pt);
  \fill (6.3, -3.8) circle (1.5pt);
  \node at (7.2,-3.5) {$\dots$};

  \node at (1.5,-0.5) {$\dots$};
  \node at (1.5,-1.5) {$\dots$};
  \node at (1.5,-3.5) {$\dots$};

 \end{tikzpicture}
\caption{
$\Ga=r\C_{su}^{(t)} +u\K_{sr,t}$, where $[{B_i}]=\C_{su}^{(t)}$ and $[{C}_j\cup {D}_j]=\K_{sr,t}$
}\label{fig-ck}
\end{figure}
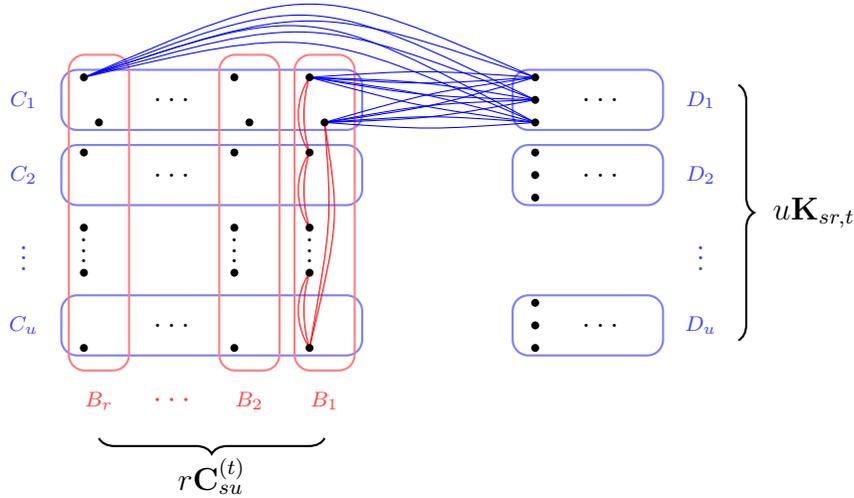
\end{center}

Each of the final family of graphs admitting a cyclic edge-bi-regular action involves three vertex-orbits. They are the graphs in the last line of Table~\ref{t:mainbireg2}, 
and Figure~\ref{fig-kk} gives a broad description of their structure, see Remark~\ref{r:figkk}. 

\begin{definition}\label{d:kk}
{\rm
Let $r,r',s,t,t',u$ be positive integers such that  
\[
\gcd(r,r')=\gcd(t,t')=\gcd(sr,ut)=\gcd(sr',ut')=1,
\]  
and let $V_1=\ZZ_{rut'}, V_2=\ZZ_{srr'}$ and $V_3=\ZZ_{r'ut}$. Define a graph $\Ga=\K\K(r,r',s,t,t',u)=(V,E)$ with 
vertex-set $V= V_1\cup V_2\cup V_3$, and edge-set $E=E_0\cup E_1$ where
\[
\begin{array}{lll}
E_0 &= \{ [i+r\ell, e_{i, \ell, k}^j, i+rk] &\mid 0\leq i\leq r-1, 0\leq \ell\leq ut'-1, 0\leq k\leq sr'-1,\\
	&			  							&\quad 0\leq j\leq t-1\},\\
E_1 &= \{ [i'+r'\ell', e_{i',\ell',k'}^{j'}, i'+r'k']&\mid 0\leq i'\leq r'-1, 0\leq \ell'\leq ut-1, 0\leq k'\leq sr-1,\\
	&										&\quad  0\leq j'\leq t'-1\} .
\end{array}
\]
Note that the last entry in each edge lies in $V_2=\ZZ_{srr'}$, while the first edge-entries lie in $V_1=\ZZ_{rut'}$ for edges in $E_0$, and in $V_3=\ZZ_{r'ut}$ for edges in $E_1$.
For $c=1,2,3$, each $v_c\in V_i$ can be written uniquely as follows:
\[
\begin{array}{lllll}
v_1	&= i+r\ell	&\mbox{\ with \ }    &0\leq i\leq r-1,	&0\leq\ell\leq ut'-1  \\ 
v_2	&= i+rk		&\mbox{\ with \ }    &0\leq i\leq r-1,	&0\leq k\leq sr'-1  \\ 
	&= i'+r'k'	&\mbox{\ with \ }    &0\leq i'\leq r'-1,&0\leq k'\leq sr-1  \\ 
v_3	&= i'+r'\ell'&\mbox{\ with \ }   &0\leq i'\leq r'-1,&0\leq\ell'\leq ut-1  \\ 
\end{array}
\]
We will define a map  $g:V\cup E \to V\cup E$ so that the vertex action is given by $v_c\to v_c+1$ where we evaluate the vertex image modulo $rut', srr', r'ut$ for $c=1,2,3$ respectively. 
To define the $g$-action on edges, set $x:=rr'sut'$ and $x':=rr'sut$, so that $|E_0|=xt$ and $|E_1|=x't'$. Then for $e_{i,\ell,k}^{j}\in E_0$ with $i,\ell,k,j$ in the ranges above, we have 
\[
1\leq (i+r\ell+1)(k+1)(j+1)\leq (rut')(sr')t=xt=|E_0|,
\] 
and similarly, for $e_{i',\ell',k'}^{j'}\in E_1$, we have 
\[
1\leq (i'+r'\ell'+1)(k'+1)(j'+1)\leq (r'ut)(sr)t'=x't'=|E_1|.
\] 
Then, writing $(i+r\ell+1)(k+1)(j+1) = ax+ b$ and $(i'+r'\ell'+1)(k'+1)(j'+1) = a'x'+ b'$, with $0\leq a\leq t-1$, $1\leq b\leq x$, $0\leq a'\leq t'-1$ and $1\leq b'\leq x'$,
we define the $g$-action by
\[
\begin{array}{lclcllll}
g:	&e_{i,\ell,k}^{j}&\to &e_{I,L,K}^{J}& \mbox{\ and \ }    &e_{i',\ell',k'}^{j'}&\to &e_{I',L',K'}^{J'}  \\ 
\end{array}
\]
where $I, J, K, L$ and $I', J', K', L'$ are as in Table~\ref{t:kk}.
}
\end{definition}

\begin{center}
\begin{table}
\begin{tabular}{cc|cccc}
$a$ 		& $i$ 			& $J$ 	& $I$ 			&$L$ 	& $K$ 	  \\  \hline 
$\leq x-2$ 	& $\leq r-2$ 	& $j$ 	& $i+1$ 		&$\ell$ & $k$ 	\\
$\leq x-2$ 	& $r-1$ & $j$ 	& $i+1$ 		&$\ell+1$ & $k+1$ 	\\ 
$x-1$ 	& $\leq r-2$ 	& $j+1$ 	& $i+1$ 		&$\ell$ & $k$ 	\\ 
$x-1$ 	& $r-1$ 		& $j+1$ 	& $i+1$ 		&$\ell+1$ & $k+1$ 	\\ \hline 
 $a'$ 		& $i'$ 		& $J'$ 	& $I'$ 		& $L'$ 		& $K'$    \\  \hline 
 $\leq x'-2$ 	& $\leq r'-2$ & $j'$ 	& $i'+1$ 	&$\ell'$ 	& $k'$    \\  
 $\leq x'-2$ 	& $r'-1$ & $j'$ 	& $i'+1$ 	&$\ell'+1$ 	& $k'+1$    \\ 
$x'-1$ 	& $\leq r'-2$ & $j'+1$ 	& $i'+1$ 	&$\ell'$ 	& $k'$    \\ 
$x'-1$ 	& $r'-1$ & $j'+1$ 	& $i'+1$ 	&$\ell'+1$ 	& $k'+1$    \\  \hline 
\end{tabular} 
\caption{Edge action for $g$ in Definition~\ref{d:kk}}\label{t:kk}
\end{table}
\end{center}

\begin{lemma}\label{l:kk}
With the notation of Definition~$\ref{d:kk}$,  the graph $\Gamma=\K\K(r,r',s,t,t',u)$ is connected, and  $\Gamma= [E_0]+[E_1]$, where $[E_0] = r\K_{sr',ut'}^{(t)}$ and $[E_1]=  r'\K_{sr,ut}^{(t')}$  as in the last line of Table~\ref{t:mainbireg2}. Further the map $g\in\Aut\Gamma$, and $\l g\r$ is cyclic of order $srr'utt'$,  with three vertex-orbits $V_1, V_2, V_3$, and is bi-regular on $E$ with orbits $E_0, E_1$.
\end{lemma}

\proof
By definition $\Gamma= [E_0]+[E_1]$ and $[E_0] = r\K_{sr',ut'}^{(t)}$ and $[E_1]=  r'\K_{sr,ut}^{(t')}$, so the first assertion holds. Further, it follows from the definition of $g$ that both  $g|_V$ and  $g|_E$ are bijections, and $g$ preserves incidence, so $g\in\Aut\Gamma$. Consider the action of $g$ on $E$. First let $e=e_{0,0,0}^0\in E_0$, and suppose that $m$ is the least positive integer such that $e^{g^m}=e$. 
It follows from Table~\ref{t:kk}, that after $xt$ repeated applications of $g$, we have $e^{g^{xt}}=e^t_{xt, y,y}$ where $y=sr'ut'$, and as $y$ is divisible by both $ut'$ and $sr'$, this edge-image is equal to $e$. Thus $g^{xt}$ fixes $e$ and so $m$ divides $xt$. Suppose for a contradiction that $0<m<xt$.  Then $m$ is of the form $m=cx+d$ where $0\leq c\leq t-1$ and $0\leq d\leq x-1$, and $e=e^{g^{m}}=e^c_{m, v,w}$, for some $v,w$. This implies that $c$ is a multiple of $t$ and hence $c=0$, so $0<m<x$. Now, with $y= sr'ut'$ so $x=ry$, we write $m=c'r+d'$ with $0\leq c'\leq y-1$ and $0\leq d'\leq r-1$. Then $e=e^{g^{m}}=e^0_{d', c',c'}$, and this implies that $d'=0$ and that $c'$ is divisible by both $ut'$ and $sr'$. Since $\gcd(sr',ut')=1$, we deduce that $c'$ is divisible by $\lcm\{sr',ut'\}= sr'ut'=y$, and hence $c'=0$ so $m=0$, a contradiction. Thus $m=xt$ and $\l g\r$ acts regularly on $E_0$. 

An analogous argument for the edge  $e^0_{0,0,0}\in E_1$ shows that $\l g\r$ also acts regularly on $E_1$, and as $|E_0|=|E_1|=xt=x't'$, the subgroup $\l g^{xt}\r$ fixes  $V\cup E$ pointwise, and hence is trivial.  Thus $\l g\r$ is edge-bi-regular, and the proof is complete.
\qed

\begin{remark}\label{r:figkk}
Figure~$\ref{fig-kk}$ gives a description of the structure of the graph 
\[
\Gamma=\K\K(r,r',s,t,t',u)= r\K_{sr',ut'}^{(t)}+r'\K_{sr,ut}^{(t')}
\] 
from Definition~\ref{d:kk}.  
The vertex-subset $V_2$ admits two $\l g\r$-invariant partitions: first the partition $\{B_1,\dots,B_{r}\}$ where the $B_i$ are the intersections with $V_2$ of the components of $[E_0]=  r\K_{sr',ut'}^{(t)}$, and second 
 $\{C_1,\dots, C_{r'}\}$ where the $C_i$ are the intersections with $V_2$ of the components of $[E_1]=r'\K_{sr,ut}^{(t')}$. The vertex-subsets $V_1, V_3$ admit the  $\l g\r$-invariant partitions
 $\{U_1,\dots, U_{r}\}$ and $\{W_1,\dots, W_{r'}\}$,  where the $U_i$ and $W_i$ are the intersections with $V_1, V_3$ of the components of $[E_0], [E_1]$, respectively.
\end{remark}

\begin{center}
  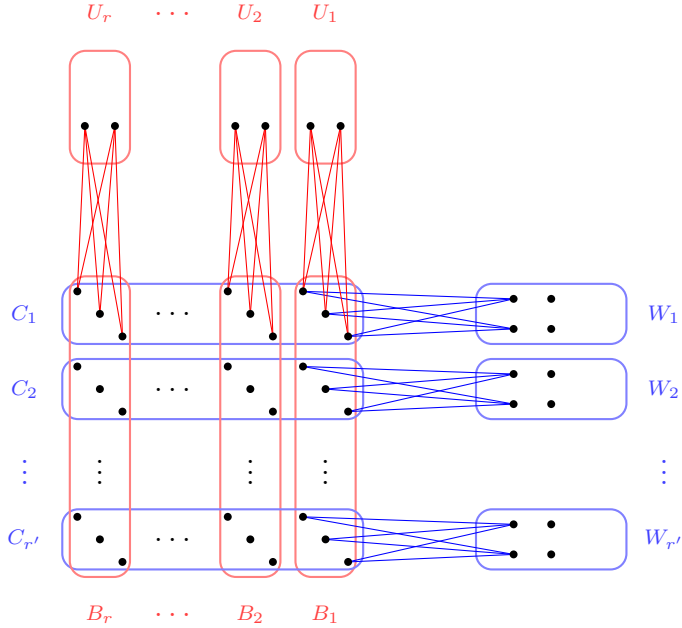
\begin{figure}[ht]
    \begin{tikzpicture}
      \foreach \x in {1,2,4} {
        \foreach \y in {1,2,4} {
          \foreach \z in {1,2,3} {
            \coordinate (v\x\y\z) at ($(-\x,-\y)+\z*(0.3,-0.3)-(0.1,-1.1) $);
          };
        };
        \foreach \y in {1,2} {
          \coordinate (w\x\y) at ($(2,-\x)+(0,1.1)-(0,\y*0.4)$);
          \coordinate (u\x\y) at ($(-\x,0) +(-0.1,0)+(\y*0.4,2)$);
        };
        \foreach \y in {3,4} {
          \coordinate (w\x\y) at ($(2.5,-\x)+(0,1.9)-(0,\y*0.4)$);
        };
      };

      \foreach \x in {1,2,4} {
        \draw[rounded corners=6pt, color= blue!50,thick] ($(0,0.9)-(0,\x)$) rectangle ($(-4,0.1)-(0,\x)$);  
        \draw[rounded corners=6pt, color= blue!50,thick] ($(1.5,0.9)-(0,\x)$) rectangle ($(3.5,0.1)-(0,\x)$);  
        \draw[rounded corners=6pt, color= red!50,thick] ($(0.9,0)-(\x,0)$) rectangle ($(0.1,-4)-(\x,0)$);  
        \draw[rounded corners=6pt, color= red!50,thick] ($(0.9,1.5)-(\x,0)$) rectangle ($(0.1,3)-(\x,0)$); 
        \foreach \y in {1,2,3} {
          \foreach \z in {1,2} {
            \draw[color=blue] (v1\x\y)--(w\x\z);
            \draw[color=red] (v\x1\y)--(u\x\z);
          };
        }; 

      };

      \foreach \x in {1,2,4} {
        \foreach \y in {1,2,4} {
          \foreach \z in {1,2,3} {
            \fill (v\x\y\z) circle (1.5pt);
          };
        };
        \foreach \y in {1,2} {
          \fill (u\x\y) circle (1.5pt);
        };
        \foreach \y in {1,2,3,4} {
          \fill (w\x\y) circle (1.5pt);
        };
        \node at ($(-2.5,.5)-(0,\x)$) {$\dots$};
        \node at ($(.5,-2.5)-(\x,0)$) {$\vdots$};
      };

      \node[color=blue!75] at (4,-.5) {\tiny $W_1$};
      \node[color=blue!75] at (4,-1.5) {\tiny $W_2$};
      \node[color=blue!75] at (4,-2.5) { $\vdots$};
      \node[color=blue!75] at (4,-3.5) {\tiny $W_{r'}$};

      \node[color=blue!75] at (-4.5,-.5) {\tiny ${C}_1$};
      \node[color=blue!75] at (-4.5,-1.5) {\tiny ${C}_2$};
      \node[color=blue!75] at (-4.5,-2.5) { $\vdots$};
      \node[color=blue!75] at (-4.5,-3.5) {\tiny ${C}_{r'}$};

      \node[color=red!75] at (-0.5,-4.5) {\tiny ${B}_1$};
      \node[color=red!75] at (-1.5,-4.5) {\tiny ${B}_2$};
      \node[color=red!75] at (-2.5,-4.5) { $\dots$};
      \node[color=red!75] at (-3.5,-4.5) {\tiny ${B}_{r}$};

      \node[color=red!75] at (-0.5,3.5) {\tiny $U_1$};
      \node[color=red!75] at (-1.5,3.5) {\tiny $U_2$};
      \node[color=red!75] at (-2.5,3.5) { $\dots$};
      \node[color=red!75] at (-3.5,3.5) {\tiny $U_{r}$};

    \end{tikzpicture}
    \caption{
$\Ga=r\K_{sr',ut'}^{(t)} +r'\K_{sr,ut}^{(t')}$, where $[{B}_i\cup U_i]=\K_{sr',ut'}^{(t)}$ and $[{C}_j\cup W_j]=\K_{sr,ut}^{(t')}$.
}\label{fig-kk}
  \end{figure}
\end{center}

\section{Coset graphs and cyclic edge-regular actions}\label{s:coset}

In \cite[Section 1.1]{LPS-rotary}, a new coset graph construction was given for arc-transitive, not necessarily simple, graphs, 
which proved helpful for analysing such graphs. Here we develop a similar theory of coset graph representations for edge-transitive bipartite graphs with two vertex-orbits. It extends the theory in \cite[Section 2]{DX} and \cite[Section 3.2]{GLP}  for simple graphs of this kind.  We use this theory in Subsection~\ref{s:cyclicedgetra} to prove Theorem~\ref{main-thm} in the case of edge-regular actions.

\begin{construction}\label{c:coset}
{\rm
Let $G$ be a group with proper subgroups $L, R, J$, such that $L\ne R$, $J\leq L\cap R$ and $J$ is 
core-free in $G$, and let $\lambda:=|L\cap R:J|$.
Define an incidence structure $\Bicos(G,L,R,J)=(V,E,\bfI)$, called a \emph{bi-coset graph}, by setting
\begin{align*}
V &= [G:L]\cup [G:R] = \{ Lx\mid x\in G\}\cup \{ Rx\mid x\in G\}\\
E &= [G:J] = \{ Jx\mid x\in G\}\\
\bfI &= \{ (Lx,Jy) \mid Lx\cap Jy\ne \emptyset\}\cup \{ (Rx,Jy) \mid Rx\cap Jy\ne \emptyset\} \subseteq V\times E.
\end{align*}
Also set $s:=|L:L\cap R|$ and $t:=|R:L\cap R|$.
}
\end{construction}

We prove that $\Bicos(G,L,R,J)$ is a $G$-edge-transitive bipartite graph and  obtain various of its other properties.

\begin{proposition}\label{p:coset}
Let $G, L, R, J, \lambda, s, t$ be as in Construction~$\ref{c:coset}$, and denote the graph $\Bicos(G,L, R, J)$ by $\Gamma$. 

\begin{enumerate}
\item[(a)] For $x, y, z\in G$, $Jy=[Lx,Jy,Rz]$ is an edge if and only if $Jy=[Ly, Jy, Ry]$, (that is, $Lx=Ly$ and $Rz=Ry$). Thus the edges of $\Gamma$ are precisely $Jy=[Ly,Jy,Ry]$, for $y\in G$.


\item[(b)] Then $\Gamma$ is a bipartite graph which is bi-regular with valencies 
$k, \ell$, and constant edge-multiplicity $\lambda$.
 
 \item[(c)] The group $G$, acting by right-multiplication, is an edge-transitive group of automorphisms of $\Gamma$ with vertex-orbits $[G:L]$ and $[G:R]$, and $L, R, J$ are the stabilisers of the vertices $L, R$ and edge $J$, respectively.
 
 \item[(d)] The base graph of $\Gamma$ is $\Gamma_0:=\Bicos(G, L, R, L\cap R)$, and $\Gamma=\Gamma_0^{(\lambda)}$. In particular $\Gamma$ is simple if and only if $L\cap R=J$.

\item[(e)]
The graph $\Gamma$ is connected if and only if $G=\langle L, R\rangle$; and $\Gamma\cong K_{s,t}^{(\lambda)}$ if and only if $G=LR$.
\end{enumerate}
\end{proposition}

\proof
(a) Let $x, y, z\in G$. Then $Jy=[Lx,Jy,Rz]$ is an edge if and only if $(Lx,Jy\in\bfI$ and $(Rz,Jy)
\in\bfI$, that is,  $Lx\cap Jy \ne \emptyset$ and $Rz\cap Jy\ne \emptyset$. This in turn is equivalent to $Lxy^{-1}\cap J \ne \emptyset$ and $Rzy^{-1}\cap J\ne \emptyset$, and since $J\leq L\cap R$, this is equivalent to $xy^{-1}\in L$ and $zy^{-1}\in R$, that is to say $Lx=Ly$ and $Rz=Ry$ so   $Jy=[Ly,Jy,Ry]$. 

(b) and (c).\ It is clear from the definition of $\Gamma$ that $[G:L]$ and $[G:R]$ form the parts of a bipartition of $\Gamma$ and each edge is incident with one vertex from each of these sets. Also, the  group $G$, acting by right-multiplication, preserves the incidence relation $\bfI$ so induces a group of automorphisms of $\Gamma$. The biparts $[G:L]$ and $[G:R]$  are the two $G$-vertex-orbits, and by part (a), $G$ is transitive on edges. The subgroups $L, R, J$ are the stabilisers of vertices $L, R$ and edge $J$ respectively, and since $J$ is core-free in $G$, $G$ acts faithfully on $\Gamma$. These transitivity properties imply that $\Gamma$ is biregular with each vertex of $[G:L]$ incident with $|L:L\cap R|=s$ vertices and  each vertex of $[G:R]$ incident with $|R:L\cap R|=t$ vertices, and that $\Gamma$ has constant edge-multiplicity $\lambda=|L\cap R:J|$ (the number of edges incident with $L$ and $R$).

(d) From the discussion in the previous paragraph it is clear that the base graph $\Gamma_0$ of $\Gamma$ has edge set identified with $[G:L\cap R]$ and is precisely the graph $\Bicos(G, L, R, L\cap R)$, and that $\Gamma$ is a $\lambda$-extender of $\Gamma_0$. 

(e) Now $\Gamma$ is connected if and only if $\Gamma_0$ is connected, and the latter is true if and only if $G=\langle L, R\rangle$, by \cite[Lemma 3.7(1)]{GLP}. Assume now that $G=LR$. Then as $R$ is the stabiliser of the vertex $R\in [G:R]$ it follows that $L$ is transitive on $[G:R]$ and similarly $R$ is transitive on $[G:L]$. Thus $L$, $R$ is adjacent to each vertex of $[G:R], [G:L]$, respectively and it follows that $\Gamma_0\cong K_{s,t}$ and $\Gamma\cong K_{s,t}^{(\lambda)}$. Conversely suppose that 
$\Gamma\cong K_{s,t}^{(\lambda)}$. Then  $\Gamma_0\cong K_{s,t}$, and as $G$ is edge-transitive, the vertex-stabiliser $L$ is transitive on the set $[G:R]$ of vertices adjacent to $L$.
Then, since $R$ is the stabiliser of $R\in[G:R]$ this means that $G=RL=LR$.
\qed

Next we show that essentially all edge-transitive graphs with two vertex-orbits arise from Construction~\ref{c:coset}.

\begin{proposition}\label{p:ceset2}
Let $\Gamma=(V, E, \bfI)$ be a graph and $G\leq \Aut\Gamma$ such that $G$ is transitive on $E$ and has two orbits $V_1, V_2$ in $V$, and $\Gamma$ has no isolated vertices, and constant edge-multuplicity $\lambda$.
Then either
\begin{enumerate}
\item[(a)] $\Gamma\cong \Bicos(G, L, R, J)$ where for some edge $e=[\a,e,\b]$ of $\Gamma$, the subgroups $L, R, J$ are the stabilisers of $\a, \b, e$ respectively, and $|L\cap R:J|=\lambda$; or

\item[(b)] $\Gamma\cong rK_2^{(\lambda)}$, where $r=|V_1|=|V_2|$.
\end{enumerate} 
\end{proposition}

\proof
If some edge $e$ of $\Gamma$ is incident with two distinct vertices of $V_i$, for some $i$, then this is true for all edges since $\Gamma$ is $G$-edge transitive and $G$ leaves $V_i$ invariant. This is a contradiction since $\Gamma$ has no isolated vertices. Hence each edge is incident with one vertex from $V_1$ and one from $V_2$. Let $e=[\a,e,\b]$ be an edge with $\a\in V_1$ and $\b\in V_2$, and let $L=G_\a, R=G_\b$ and $J=G_e$. Then we may identify the sets $V_1, V_2, E$ with $[G:L], [G:R], [G:J]$, respectively, with $G$ acting by right multiplication. Note that $G_e$ fixes the unique vertex of $V_i$ with which it is incident, for each $i$, and so $J=G_e\leq G_\a\cap G_\b=L\cap R$. With this identification, the edge $e=[L, J, R]$ and each edge can be expressed as $e^y=[Ly, Jy, Ry]$ for some $y\in G$ (since $G$ is transitive on $E$). Thus $Lx$ is incident with $Jz$ if and only if there exists $y$ such that $Lx=Ly$ and  $Jz=Jy$, and hence $Lx\cap Jz$ contains $y$ so is non-empty. Conversely, if $Lx\cap Jz$ contains an element $y$, then $Lx=Ly$ and  $Jz=Jy$. Similar statements hold for incidences between $Jx$ and $Rz$. Thus, provided $L\ne R$, we see that, under these identifications we have $\Gamma= \Bicos(G, L, R, J)$ and part (a) holds. 

Finally suppose that $L=R$. Since $G$ is transitive on $E$, $L$ is transitive on the edges incident with $\a$, and since $L=R$ fixes $\b$ this means that all of these edges are incident with $\b$. Thus the connected component of $\Gamma$ induced on $\{\a,\b\}$ is $K_2^{(\lambda)}$, and we have $\Gamma\cong rK_2^{(\lambda)}$, where $r=|V_1|=|V_2|$, as in (b).
\qed

\subsection{Graphs admitting cyclic edge-transitive groups}\label{s:cyclicedgetra}

We use the theory in the previous subsection, to determine all graphs admitting a cyclic subgroup of automorphisms that is regular on edges. First we consider connected graphs. We assume that $|V|\geq3$, in the light of Proposition~\ref{p:k2}.

\begin{theorem}\label{th:cyclicedgetra}
Let $\Gamma=(V, E, \bfI)$ be a connected graph with $|V|\geq3$, let $G\leq \Aut\Gamma$ be a cyclic subgroup acting regularly on $E$, and let $\Gamma$ have edge-multiplicity $\lambda$. Then $\Gamma$ and $|G|$ are as in one of the lines of Table~$\ref{t:cyclicedgetra}$, and the induced permutation groups $G^X$ on vertices ($X=V$), edges ($X=E$) and arcs ($X=A$) are transitive if the entry in the column headed $G^X$ is $\surd$ and intransitive if the entry is $\times$. 
\end{theorem}

\begin{center}
\begin{table}
\begin{tabular}{cccccl}
$\Gamma$ & $|G|$ & $G^V$ & $G^E$ & $G^A$ & Conditions \\ 
\hline 
$\C_n^{(\lambda)}$ & $n\lambda$ & $\surd$ & $\surd$ & $\times$ & $n\geq 3$ \\ 
$\K_{s,t}^{(\lambda)}$ & $st\lambda$ & $\times$ & $\surd$ & $\times$ & $\gcd(s,t)=1$, $st>1$ \\ 
\hline 
\end{tabular} 
\caption{Table for Theorem~\ref{th:cyclicedgetra}}\label{t:cyclicedgetra}
\end{table}
\end{center}

\proof
Since $|V|\geq3$ we have $\Gamma\not\cong\K_2^{(\lambda)}$, and since $\Gamma$ is connected and $G$-edge-transitive, $G$ has at most two orbits in $V$. Suppose first that $G$ has two vertex-orbits, $V_1$ and $V_2$.  Then by Proposition~\ref{p:ceset2}, 
$\Gamma=\Bicos(G,L,R,J)$ with $L=G_\a, R=G_\b, J=G_e$ for some edge $e=[\a,e,\b]$. 
It follows from Proposition~\ref{p:coset}(e) that $G=\langle L, R\rangle$, and since $G$ is cyclic this means that $G=LR$ and hence that $\Gamma=\K_{k,\ell}^{(\lambda)}$, where $s=|L:L\cap R|$ and $t=|R:L\cap R|$,
and  $st>1$ since $\Gamma\not\cong \K_2^{(\lambda)}$. Also $L\cap R$ acts trivially on $V$ and the edge-stabiliser $J$ acts trivially on $E$ (since $G$ is abelian), and hence $J=1$ by Lemma~\ref{l:aut}. Thus $|G|=|E|=st\lambda$,
and $|L\cap R|=\lambda$. Now we have $|L|=s\lambda$ and $|R|=t\lambda$ and so $|G|=
|LR|=\lcm\{s\lambda, t\lambda\}=\lambda\,\lcm\{s, t\}$, and hence $st=\lcm\{s,t\}$, that is, $\gcd(s,t)=1$. Thus all the entries of line 2 of Table~\ref{t:cyclicedgetra} hold.  

We may now assume that $G$ is transitive on $V$ as well as on $E$. Suppose that $G$ is not transitive on arcs. Then the stabiliser of an edge $e=[\a,e,\b]$ satisfies $G_e=1$, and $|G|=|E|$. Also, the stabiliser $G_\a$ is normal in $G$ since $G$ is abelian, and so $G_\a$ acts trivially on $V$. This implies that $G_\a=G_\b$ and so $G_\a$ is transitive on the $\lambda$ edges incident with $\a$ and $\b$. Since $G_e=1$ and $G_e\leq G_{\a,\b}=G_\a$,  we have $|G_\a|=\lambda$, and 
$n:=|V|=|G:G_\a|$ so $|G|=n\lambda$, and note that $n=|V|\geq3$.  Let $\Gamma_0$ be the base graph of $\Gamma$ so $\Gamma=\Gamma_0^{(\lambda)}$. Then $G/G_\a$ acts transitively on the vertices and edges of $\Gamma_0$, and  $\Gamma_0$ has $|E|/\lambda = n = |V|$ edges and also $n$ vertices. These properties imply that $\Gamma_0$ has valency $2$ and hence is a cycle of length $n$. Hence $\Gamma=\C_n^{(\lambda)}$ and $n\geq 3$, so all the entries of line 1 of Table~\ref{t:cyclicedgetra} hold. 

Suppose finally that $G$ is transitive (hence regular) on the arcs of $\Gamma$. Then the stabiliser $G_e$ of an edge $e=[\a,e,\b]$ contains an element $g$ which interchanges $\a$ and $\b$. Now $G_e\cap G_\a$ acts trivially on both $V$ and $E$ (since $G$ is abelian and transitive on $V$ and $E$), and hence $G_e\cap G_\a=1$. Thus $G_e=\langle g\rangle \cong Z_2$. Moreover, $g\not\in G_\a$, but $g$ normalises $G_\a$ (since $G$ is abelian) and $G=\langle G_\a,g\rangle$ (since $\Gamma$ is connected). Thus $|V|=|G:G_\a|=2$, which is a contradiction.
\qed

Now we use this classification to determine all  graphs with  no isolated vertices which admit a cyclic subgroup of automorphisms regular on edges. Such graphs have constant 
edge-multiplicity.

\begin{corollary}\label{cyc-edge-reg-2}
Let $\Gamma=(V, E, \bfI)$ be a graph with no isolated vertices, and constant edge-multiplicity $\lambda$, and let $G\leq \Aut\Gamma$ be a cyclic subgroup transitive 
on $E$. Then  
 \begin{itemize}
 \item[(a)] $\Gamma = r\Gamma_0$ where $\Gamma$ has $r$ connected components, each isomorphic to $\Gamma_0$, and
 \item[(b)] for $G_0$ the setwise stabiliser in $G$ of a component $\Gamma_0$, $|G|=r |G_0|$, and either $\Gamma_0=\K_2^{(\lambda)}$, or $(\Gamma_0, |G_0|)$ is as in one of the lines of Table~\ref{t:cyclicedgetra}, and $\Ga_0$ is as in one of the lines of Table~\ref{t:mainreg}.
 \end{itemize}
\end{corollary}

\proof
The result follows from Theorem~\ref{th:cyclicedgetra} if $\Ga$ is connected, so suppose that $\Ga$ has $r\geq 2$ connected components $\Sigma_1,\dots, \Sigma_r$, and let $E_i$ be the edge-set of $\Sigma_i$, for $1\leqslant i\leqslant r$.
Since $G$ is edge-transitive it follows that 
$\Gamma_0:=\Sigma_1\cong \dots\cong \Sigma_r$, so $\Gamma = r\Gamma_0$,    and $\{E_1,\dots,E_r\}$ forms a 
$G$-invariant partition of $E$.
Then since $G$ is cyclic, it induces a regular action on $\{E_1,\dots,E_r\}$, and the unique subgroup $G_0$ of index $r$ in $G$ is the setwise stabiliser of $E_i$ for each $i$.
Also $G_0$ is transitive on each $E_i$ and $G_0$ is cyclic, and hence the unique subgroup $J$ of $G_0$ of index $|E_i|$ is the stabiliser of each edge of each of the $E_i$. That is to say, $J$ is the kernel of the action of $G$ on $E$. It follows from Lemma~\ref{l:aut} that either  $J=1$, or  $\Gamma_0=\K_2^{(\lambda)}$. Thus we may assume that $J=1$. In particular $G_0$ acts edge-regularly and faithfully as a cyclic group of automorphisms of $\Sigma_i$, for each $i$. Then by Theorem~\ref{th:cyclicedgetra}, $(\Gamma_0, |G_0|)$ is as in one of the lines of Table~\ref{t:cyclicedgetra}, and so $\Ga_0$ is as in one of the lines of Table~\ref{t:mainreg}.
\qed

\section{Proof of Theorem~\ref{main-thm}}

In this section we complete the proof of Theorem~\ref{main-thm}. 
Let $\Ga=(V,E)$ be a connected graph with $|V|\geq3$, and assume that  $G\leqslant\Aut\Ga$ is a cyclic subgroup which is regular or bi-regular on the edge set $E$. In particular $\Ga$ has no isolated vertices. If $G$ is regular on $E$ then the possibilities for $(\Gamma, G)$ are determined in Theorem~\ref{th:cyclicedgetra} and (see Corollary~\ref{cyc-edge-reg-2}(b)) $\Ga$ is as in one of the lines of Table~\ref{t:mainreg}.
So we may assume that  $G$ is bi-regular on $E$, with two edge orbits $E_0,E_1$ of equal size. Let
\[
\Pi_i=[E_i] = (V(E_i),E_i),\ \mbox{where $i=0$ or 1}
\]
be the induced subgraphs.
Then $\Pi_0,\Pi_1$ are edge disjoint graphs, and by our convention, $\Gamma=\Pi_0 + \Pi_1$. By the definition of an induced subgraph, $\Pi_i$ has no isolated vertices, and by the definition of bi-regular, $G$ induces an edge-regular action on $\Pi_i$ for each $i$. Hence
by Corollary~\ref{cyc-edge-reg-2}, the following holds for each of $i=0$ or 1. The graph $\Pi_i=r_i\Sigma_i$, where $\Pi_i$ has $r_i$ connected components, each isomorphic to $\Sigma_i$, and for $H_i$ the setwise stabiliser in $G$ of a component $\Sigma_i$, $|G|=r_i|H_i|$, and one of 

\begin{itemize}
\item[(1)] $\Sigma_i=\C_{n_i}^{(\lambda_i)}$, and $H_i=\ZZ_{n_i\lambda_i}$ is transitive on $V(E_i)$, for some $n_i\geqslant3, \lambda_i\geq1$; or
\item[(2)] $\Sigma_i=\K_{s_i,t_i}^{(\lambda_i)}$, and $H=\ZZ_{s_it_i\lambda_i}$ is bi-transitive on $V(E_i)$, with $s_it_i>1$ and $\gcd(s_i,t_i)=1, \lambda_i\geq1$; or
\item[(3)] $\Sigma_i=\K_{2}^{(\lambda_i)}$, and by Proposition~\ref{p:k2}, one of lines 1--3 of Table~\ref{t:k2} holds for $H_i$.
\end{itemize}

It follows from the above discussion (and Corollary~\ref{cyc-edge-reg-2}) that $G$ has at most two orbits in $V(E_i)$ for each $i$, and since $\Gamma$ is connected, the sets $V(E_0)$ and $V(E_1)$ are not disjoint. Thus $V(E_0)$ and $V(E_1)$ share at least one $G$-orbit, and hence $G^V$ has at most three orbits. We show that $\Gamma$ and $|G|$ are as in one of the cases of Theorem~\ref{main-thm} in the following three subsections, according to the number of $G$-orbits in $V$.

\medskip

\subsection{The case where $G^V$ is transitive}\label{s:T1-1-vert-trans}

Suppose in this subsection that $G^V$ is transitive, so $V=V(E_i)$ for each $i$.

Assume first that $\Pi_0= r_0 \K_2^{(\lambda_0)}$. Then, by Proposition~\ref{p:k2}, line 1 or 2 of Table~\ref{t:k2} holes, so either (i) $|G|=2|E_0|=2r_0\lambda_0$  with $\lambda_0$ odd, and the $G$-action on $\Pi_0$ is arc-transitive, or (ii)  $|G|= |E_0|=r_0\lambda_0$ with $\lambda_0$ even and $G$ is faithful on $E_0$ and not arc-transitive on $\Pi_0$. In either case $|G|$ is even, so $G$ has a unique subgroup $X$ of order $2$, and the $X$-orbits in $V$ are the components of $\Pi_0$. Suppose to start with that $\Pi_1= r_1 \K_2^{(\lambda_1)}$. Then $|V|=2r_0=2r_1$ so $r_0=r_1$, and $|E_i|=r_0\lambda_0=r_1\lambda_1$ (by definition, since $G$ is biregular on $E$), so also $\lambda_0=\lambda_1=\lambda$, say, and $|G|=2r_1\lambda$ or $r_1\lambda$ according as 
$\lambda$ is odd or even, respectively. It follows that the $X$ vertex-orbits are also the connected components of $\Pi_1$, so a connected component of $\Gamma$ has size $2$. Since $\Gamma$ is connected this implies that $|V|=2$, which is a contradiction. Hence $\Pi_1\neq  r_1 \K_2^{(\lambda_1)}$. This implies in particular that 
$G$ acts faithfully on $E_1$ (by Lemma~\ref{l:aut}). Thus $|G|=|E_1|$ since $G$ is abelian, and therefore, since $|E_0|=|E_1|$ as $G$ is bi-regular on $E$, we have $|G|=|E_0|$ so (ii) above holds, that is, $|G|=r_0\lambda_0$, $\lambda_0$ is even, and $G$ is faithful on $E_0$.  
Also, since $\Pi_1\ne r_1\K_2^{(\lambda_1)}$ and $G^V$ is transitive, it follows from Corollary~\ref{cyc-edge-reg-2} that  
$\Pi_1= r_1 \C_{n_1}^{(\lambda_1)}$ for some $n_1\geq3$. Then $|V|=2r_0=r_1n_1$, and $|E_i|=r_0\lambda_0=r_1n_1\lambda_1$, so $\lambda_0=2\lambda_1$. 

Suppose that $X$ fixes a $\Pi_1$-component $\Sigma$ setwise. Then $V(\Sigma)$ is a union of (the vertex-sets of) some components of $\Pi_0$, and since $\Gamma$ is connected it follows that  
$\Pi_1=\Sigma$, so $V=V(\Sigma)$,  $r_1=1$ and $n_1=2r_0$ is even. Moreover the two vertices of each component $\K_2^{(\lambda_0)}$ of $\Pi_0$ form an antipodal pair of vertices of  $\Pi_1=\C_{n_1}^{(\lambda_1)}$, and $X$ interchanges this vertex-pair and interchanges in pairs the $\lambda_0=2\lambda_1$ edges of $\Pi_0$ incident with them.
Thus, setting $\lambda:=\lambda_1$ and $n :=n_1/2$, we have  $\Gamma= \Gamma_0^{(\lambda)}$ with $|G|=2n\lambda$ and $\Gamma_0$ the Cayley graph $\Gamma(2n, 1, n)=\C_{2n}+n\K_2^{(2)}$ of Corollary~\ref{c:circ1}(a), and line 2 of Table~\ref{t:mainbireg1} holds. 

On the other hand, suppose that $X$ interchanges the components of $\Pi_1$ in pairs. 
Then $r_1$ is even  and the two vertices of a component $\K_2^{(\lambda_0)}$ of $\Pi_0$ lie in two distinct components of $\Pi_1$. Further, the union of these two distinct components of $\Pi_1$ is also a union of $n_1$ components of $\Pi_0$ and hence forms a component of $\Gamma$. Since $\Gamma$ is connected we conclude that $r_1=2$ and $r_0=n_1=n$, say. Again the group $X$ interchanges the two $\Pi_1$-components $\C_{n}^{(\lambda)}$, where $\lambda:=\lambda_1$, and for each component $\K_2^{(2\lambda)}$ of $\Pi_0$, $X$ interchanges the $2\lambda$ edges in pairs.  Since $G$ is cyclic and transitive on the $2n=|V|$ vertices of $\Gamma$, and since the index $2$ subgroup of $G$ stabilising a component of $\Pi_1$ is transitive on the $n$ vertices of that component, it follows that $n$ is odd, and hence  $|G|=2n\lambda$ and  $\Gamma= \Gamma_0^{(\lambda)}$ with $\Gamma_0$ the Cayley graph $\Gamma(2n, 2, n)=2\C_{n}+n\K_2^{(2)}$ of Corollary~\ref{c:circ1}(b), and line 2 of Table~\ref{t:mainbireg1} holds. 
 
Thus we may assume from now on that no component of either $\Pi_0$ or  $\Pi_1$ is $\K_2^{(\lambda)}$ for any $\lambda$. Then by Corollary~\ref{cyc-edge-reg-2}, since $G^V$ is transitive, each $\Pi_i= r_i \C_{n_i}^{(\lambda_i)}$ for some $n_i\geq3$. Thus $|V|=r_1n_1=r_0n_0$ and  $|E_1|=|E_0|=r_1n_1\lambda_1=r_0n_0\lambda_0$, and hence in particular  $\lambda_1=\lambda_0$. This means that $\Gamma = \Gamma_0^{(\lambda_0)}$ and $\Pi_i=\Phi_i^{(\lambda_0)}$, for  graphs $\Phi_i$ with edge-multiplicity $1$ admitting a cyclic vertex-transitive  group $H$ (induced by $G$)  of order $n:=r_1n_1=r_0n_0$. Thus $H$ is bi-regular on the edge-set of $\Gamma_0$ with orbits the edge-sets of $\Phi_0$ and $\Phi_1$.  This means that each 
$\Phi_i$ may be regarded as a simple Cayley graph for the cyclic group $H$. If as simple graphs, the edges sets of $\Phi_0$ and $\Phi_1$ coincide then, since $\Gamma$ is connected, $r_1=r_0=1$ and $n=n_0=n_1$, and setting $\lambda:=\lambda_0$, $\Gamma = \Circ(n,\{1,-1\})^{(2\lambda)}=\C_n^{(2\lambda)},$ with $G=\ZZ_{n\lambda}$ acting with two edge orbits $E_0$ and $E_1$. Thus
$\Gamma = \Gamma_0^{(\lambda)}, \mbox{with $\Gamma_0= \C_n^{(2)}$ and $|G|=n\lambda$}$, and line 1 of Table~\ref{t:mainbireg1} holds.

If this is not the case then the edge sets of $\Phi_0$ and $\Phi_1$ are disjoint (as they are $G$-orbits),  and hence 
$\Gamma_0=\Circ(n,S)$, where $S=\{ a, -a, b, -b\}\subset \ZZ_n\setminus\{0\}$ with $\Phi_0=\Circ(n,\{ a,-a\})$ and $\Phi_1=\Circ(n,\{ b,-b\})$, for some $a, b$ with $|a|\geq 3, |b|\geq 3$, $a\ne \pm b$,  (see Section~\ref{s:circ}).    Since $\Gamma$, and hence also $\Gamma_0$ is connected it follows (as in the proof of Lemma~\ref{l:circ1}) that $\gcd(n,a,b)=1$, and hence, setting $\lambda:=\lambda_1=\lambda_0$, $\Gamma = \Gamma_0^{(\lambda)}$, where $\Gamma_0=\Circ(n, \{a,-a, b,-b\})$ and $|G|=n\lambda$  as in 
line 6 of Table~\ref{t:mainbireg1} (see Lemma~\ref{l:circ1}).

This proves Theorem~\ref{main-thm} for cyclic vertex-transitive, edge-bi-regular actions.

\subsection{The case where $G^V$ has two orbits}\label{s:T1-1-twovertorbs}

In this subsection we assume that $G$ has two orbits $V_1, V_2$  on vertices. Suppose first that $V=V(E_i)$ for each $i$. Then $G$ acts faithfully on each $E_i$. Assume that $\Pi_0= r_0 \K_2^{(\lambda_0)}$, so by Proposition~\ref{p:k2},  $|V|=2r_0$, $|G|=|E_0|=|E_1|=r_0\lambda_0$, and each of the two $G$-orbits in $V$ has size $r_0$. We claim that also $\Pi_1= r_1 \K_2^{(\lambda_1)}$. If this is not the case then, by Corollary~\ref{cyc-edge-reg-2}, $\Pi_1 = r_1 \K_{s,t}^{(\lambda_1)}$ with $\gcd(s,t)=1$ and $st>1$. However this means that the $G$-vertex-orbits have unequal sizes $r_1s$ and $r_1t$, which is a contradiction. Hence $\Pi_1= r_1 \K_2^{(\lambda_1)}$, and so $|V|=2r_1$. Thus $r_0=r_1=r$, say, and $|E_1|=r_1\lambda_1$ so $\lambda_0=\lambda_1=\lambda$, say. If some $\Pi_0$-component is also a $\Pi_1$-component, then $r=1$ since $\Gamma$ is connected, and $|V|=2$, which is a contradiction. Thus no $\Pi_0$-component is equal to a $\Pi_1$-component. Hence $\Gamma=\Gamma_0^{(\lambda)}$, and the graph $\Gamma_0$ is connected (since $\Gamma$ is connected) with edge-multiplicity $1$, and valency 2 (since each vertex is adjacent to exactly one vertex by an edge from a $\Pi_i$-component, for each $i$). Thus $\Gamma_0=\C_n$, where $n=r_0+r_1=2r$, and the group $G$ induces a cyclic subgroup $H\leq\Aut\Gamma_0$ of order $n/2=r$ with two vertex-orbits and two edge-orbits. Thus, as in Lemma~\ref{l:n-cycle}, $\Gamma =\C_n^{(\lambda)},$ with $\$n$ even, $\Pi_0\cong\Pi_1\cong (n/2)\,\K_2^{(\lambda)}$ and  $|G|=n\lambda/2$, and line 1 of Table~\ref{t:mainbireg2} holds for $\Gamma_0$.

Therefore, in the case where $V=V(E_i)$ for each $i$,  we may assume that neither of the $\Pi_i$ has a component $\K_2^{(\lambda)}$ for any $\lambda$. It follows from Corollary~\ref{cyc-edge-reg-2} that, for each $i$, $\Pi_i = r_i \K_{s_i,t_i}^{(\lambda_i)}$ with $\gcd(s_i,t_i)=1$ and $s_it_i>1$. Without loss of generality we may assume that the $G$-vertex-orbits $V_1, V_2$ have sizes $|V_1|=r_0s_0=r_1s_1$ and $|V_2|=r_0t_0=r_1t_1$. Also we have $r_0s_0t_0\lambda_0=|E_0|=|E_1|=
r_1s_1t_1\lambda_1$. Let $\Sigma_1$ be a $\Pi_1$-component. Suppose first that, for some vertices $u_i\in V_i\cap V(\Sigma_1)$, for $i=1, 2$, there exists an edge $e_0=[u_1, e_0, u_2]\in E_0$. Since $\Sigma_1 = \K_{s_1,t_1}^{(\lambda_1)}$, we also have an edge $e_1= [u_1, e_1, u_2]\in E_1$, and by the transitivity of $G$ on $E_0$ and $E_1$, it follows that $V(\Sigma_1)$ is contained in the vertex-set of a $\Pi_0$-component $\Sigma_0$, and the same argument gives $V(\Sigma_0)\subseteq V(\Sigma_1)$, yielding equality $V(\Sigma_1)=V(\Sigma_0)$. This implies that $s:=s_1=s_0, t:=t_1=t_0$, and since $\Gamma$ is connected, also $r_1=r_0=1$ and hence $\lambda_1=\lambda_0=\lambda$, say. Thus $\Gamma = \Gamma_0^{(\lambda)}$, with $\Gamma_0=\K_{s,t}^{(2)}$ and $|G|=st\lambda$, and line 2 of Table~\ref{t:mainbireg2} holds for $\Gamma_0$.

We may therefore assume that, for $u\in V_1\cap V(\Sigma_1)$, the set $\Pi_1(u)$ of $t_1$ vertices adjacent to $u$ in $\Pi_1$ (that is to say, adjacent in $\Sigma_1$) is disjoint from the set $\Pi_0(u)$ of $t_0$ vertices adjacent to $u$ in $\Pi_0$.   Since $G$ is transitive on $E_1$ and $E_0$, it follows that $G_u$ is transitive on each of the disjoint sets $\Pi_1(u)$ and $\Pi_0(u)$. Moreover since $G$ is cyclic, $G_u$ is normal in $G$ and hence all of its orbits in $V_2$ have the same size. Hence   $t_1=t_0=t$, say, and the number of $G_u$-orbits in $V_2$ is  $r_0=|V_2|/t_0=|V_2|/t_1=r_1 = r,$ say, and $r\geq2$. An analogous argument with a vertex $w\in V_2\cap V(\Sigma_1)$ yields that $s_0=s_1=s$, say, and $G_w$ has all orbits in $V_1$ of length $s$. Thus  $\lambda_0=|E_0|/rst = |E_1|/rst = \lambda_1=\lambda$, say, and it follows that $\Pi_0\cong\Pi_1\cong r\, \K_{s,t}^{(\lambda)}$. 
Since $E_0, E_1$ are disjoint, it follows that $\Gamma=\Gamma_0^{(\lambda)}$, and each $\Pi_i=\Phi_i^{(\lambda)}$, with $\Gamma_0=(V,F)$ connected and $\Phi_i=(V,F_i)\cong r\,\K_{s,t}$ such that $E=F^{(\lambda)}$, each $E_i=F_i^{(\lambda)}$, and $F$ is the disjoint union $F_0\cup F_1$. Also $G$ induces a cyclic bi-transitive, edge-bi-regular group $H$ on $\Gamma_0$. Note that the $H$-action on $V$ is equivalent to the $G$-action on $V$, and $|H|=rst$. 

Let $L=H_u$ and $K=H_w$ so the $K$-orbits in $V_1$  and the $L$-orbits in $V_2$ are the vertex-subsets, in $V_1, V_2$ respectively, of the components of the $\Pi_i$, or equivalently the $\Phi_i$. Since $H=\l h\r$ permutes these two families of $r$ subsets cyclically, we may label the subsets as $U_\ell$, with $\ell\in\ZZ_{2r}$, such that $V_1=\cup\{ U_{2k}\mid 0\leq k\leq r-1\}$ , $V_2=\cup\{ U_{2k+1}\mid 0\leq k\leq r-1\}$, with each $|U_{2k}|=s, |U_{2k+1}|=t$, and $U_{2k}\cup U_{2k+1}$ the vertex set of a $\Phi_1$-component, and such that $U_\ell^h=U_{\ell+2}$ for each $\ell\in\ZZ_{2r}$. Now the vertex-set of the $\Phi_0$-component containing $U_1$ is $U_{2u}\cup U_{1}$, for some $2u\in\ZZ_{2r}$, and the $H$-action implies that the $\Phi_0$-components have vertex-sets $U_{2(k+u)}\cup U_{2k+1}$, for $0\leq k\leq r-1$. It follows that there are paths in $\Gamma_0$ from vertices in $U_0$ (with edges alternately in $F_1$ and $F_0$) to vertices in $U_{2u\ell}$ for all $\ell$, and to no other vertices in $V_1$. Since $\Gamma_0$ is connected, this implies that $\gcd(r,u)=1$, and hence there exists $v\in\ZZ_r$ such that $uv\equiv 1\pmod{r}$.  Let $\varphi:V\to V$ be any map which induces bijections $U_{2k}\to U_{2kv}$ and $U_{2k+1}\to U_{2kv+1}$ for each $k$. Then  $\varphi$ permutes the 
$\Phi_1$-components among themselves, and maps the $\Phi_0$-component with vertex-set  $U_{2(k+u)}\cup U_{2k+1}$ to a graph $\K_{s,t}$ with vertex set  $U_{2kv+2}\cup U_{2kv+1}$. Thus  $\varphi$ induces a graph isomorphism from $\Gamma_0$ to the graph $\C_{2r}[s\K_1,t\K_1]$ in Definition~\ref{d:cyclest}. Therefore 
$\Gamma = \Gamma_0^{(\lambda)}$, with $\Gamma_0\cong \C_{2r}[s\K_1,t\K_1]$,  $|G|=rst\lambda$, and line 3 of Table~\ref{t:mainbireg2} holds for $\Gamma_0$.

This completes our analysis of the case where $V=V(E_i)$ for each $i$, so we may without loss of generality assume from now on that $V(E_1)=V_1$ and (since $\Gamma$ is connected) that $V(E_0)=V=V_1\cup V_2$. In particular $|V_1|>1$ since $E_1\ne\emptyset$. Then, by Corollary~\ref{cyc-edge-reg-2}, we have (i) $\Pi_1=r_1\,\Sigma_{1}^{(\lambda_1)}$ with $\Sigma_{1}= \K_2$ or $\C_{n}$ and $|V_1|=2r_1$ or $nr_1$ respectively, and (ii) $\Pi_0= r_0 \Sigma_0^{(\lambda_0)}$ with $\Sigma_0= \K_2$ or $\K_{s,t}$ where $\gcd(s,t)=1$ and $st>1$. For case (ii) we remove the constraint $st>1$ and consider the two possibilities for $\Sigma_0$ together, and further, we assume that each $\Pi_0$-component has $s\geq 1$ vertices in $V_1$ and $t\geq 1$ vertices in $V_2$, where $\gcd(s,t)=1$, so $|V_1|=r_0s$ and $|V_2|=r_0t$.

Suppose first that some $\Pi_1$-edge (that is, an edge in $E_1$) is incident with two distinct vertices from the same $\Pi_0$-component, so in particular  $s>1$. Since $G$ is transitive on $E_0$ and acts as automorphisms of $\Pi_1$, this holds for all edges in $E_1$, and hence there are no $E_1$-edges between distinct $\Pi_0$-components. Since $\Gamma$ is connected this implies that $\Pi_0$ is connected, so $r_0=1$, $|V_2|=t$, and $|V_1|=s=2r_1$ or $nr_1$ according as $\Sigma_{1}= \K_2$ or $\C_{n}$, respectively. Consider first $\Sigma_{1}= \K_2$, so $s=2r_1$ and hence $t$ is odd and $\gcd(t,r_1)=1$. Then $|E_1|=r_1\lambda_1$ and $|E_0|=st\lambda_0=2r_1t\lambda_0$, so with $\lambda:=\lambda_0$ and $r:=r_1$ we have $\lambda_1=2t\lambda$. Thus $\Pi_0\cong \K_{2r,t}^{(\lambda)}$ and $\Pi_1=r\,\K_2^{(2t\lambda)}$, and as in Definition~\ref{d:ck},  
$\Gamma = \Gamma_0^{(\lambda)}$ with $\Gamma_0=r\,\K_2^{(2t)}+\K_{2r,t}$, $|G|=2rt\lambda$, and line 4 of Table~\ref{t:mainbireg2} holds for $\Gamma_0$.

Now consider $\Sigma_{1}= \C_n$, so $s=nr_1$ and $\gcd(t,nr_1)=1$. Then $|E_1|=nr_1\lambda_1$ and $|E_0|=nr_1t\lambda_0$, so with $\lambda:=\lambda_0$ and $r:=r_1$ we have $\lambda_1=t\lambda$, $\Pi_0\cong \K_{nr,t}^{(\lambda)}$ and $\Pi_1=r\,\C_n^{(t\lambda)}$, and as in Definition~\ref{d:ck},  
$\Gamma = \Gamma_0^{(\lambda)}$ with $\Gamma_0=r\,\C_n^{(t)}+\K_{nr,t}$, $|G|=nrt\lambda$, 
 and line 6 of Table~\ref{t:mainbireg2} holds for $\Gamma_0$.

Finally suppose that each edge of $E_1$ is incident with vertices from two different components of $\Pi_0$, so $r_0\geq2$. Thus $e\in E_1$ satisfies $e=[\a_1,e,\a_2]$ with 
$\a_i$ in component $\Sigma_{0,i}$ of $\Pi_0$ for $i=1,2$, and the subgroup $H$ of index $r_0$ in $G$ fixes each of the $\Pi_0$-components setwise (since $H$ is normal in $G$), and the $H$-orbits in $V_1$ are the $s$-subsets of $V_1$ lying in the  $\Pi_0$-components. 
Thus each of the $s$ vertices of $\Sigma_{0,1}$ in $V_1$ is joined by an edge of $E_1$ to a vertex in $\Sigma_{0,2}$. Further since $G$ permutes the $\Pi_0$-components transitively and cyclically, it follows that the component, $\K_2^{(\lambda_1)}$ or $\C_n^{(\lambda_1)}$, of $\Pi_1$ containing $e$ meets each of the $\Pi_0$-components. Consider first $\Sigma_1=\K_2$. Then there are exactly $r_0=2$ components of $\Pi_0$, and $r_1=s$  components of $\Pi_1$. Thus $|V_1|=2s$ and for the transitive cyclic group $\ZZ_{2s}$ induced on $V_1$, the stabiliser $\ZZ_2$ of a $\Pi_1$-component interchanges its two vertices and hence interchanges the two $\Pi_0$-components. This implies that $s$ is odd. Also $s\lambda_1=|E_1|=|E_0|=2st\lambda_0$, so setting $\lambda:=\lambda_0$ we have $\lambda_1=2t\lambda$ and so 
$\Gamma = \Gamma_0^{(\lambda)}$, with $\Gamma_0=s\,\K_2^{(2t)}+2 \K_{s,t}$, $|G|=2st\lambda$,  and line 5 of Table~\ref{t:mainbireg2} holds for $\Gamma_0$.

Now consider $\Sigma_{1}= \C_n$. Then $\Sigma_1$ meets each of the $\Pi_0$-components in a constant number $u$ of points, so we have $n=ur_0$ and $s=ur_1$, so $\gcd(ur_1,t)=1$. Now the stabiliser of a $\Pi_1$-component is still transitive on the $\Pi_0$-components, and conversely, and this holds if and only if $\gcd(r_0,r_1)=1$. Also  $|E_1|=nr_1\lambda_1=ur_0r_1\lambda_1$ and $|E_0|=r_0st\lambda_0=r_0ur_1t\lambda_0$, so with $\lambda:=\lambda_0$ we have $\lambda_1=t\lambda$, and so $\Pi_1=r_1\,\C_{ur_0}^{(t\lambda)}, \Pi_0=r_0\K_{ur_1,t}^{(\lambda)}$, $|G|=ur_0r_1t\lambda$,  and
$\Gamma = \Gamma_0^{(\lambda)}$ with $\Gamma_0=r_1\,\C_{ur_0}^{(t)}+r_0\K_{ur_1,t}$, $ur_0\geq3$, $r_0\geq2$, and $\gcd(r_0,r_1)=\gcd(ur_1,t)= 1$. So line 7 of Table~\ref{t:mainbireg2} holds for $\Gamma_0$.

This completes the analysis of the case where there are two $G$-vertex-orbits.

\subsection{The case where $G^V$ has three orbits}\label{s:T1-1-threevertorbs}

In this final subsection we assume that $G$ has three orbits $V_1, V_2, V_3$  on vertices with, say, $V(E_0)=V_1\cup V_2$ and  $V(E_1)=V_2\cup V_3$. Then $G$ acts faithfully on each $E_i$, and by Corollary~\ref{cyc-edge-reg-2}, we have $\Pi_0= r_0\K_{s_0, t_0}^{(\lambda_0)}$ and $\Pi_1= r_1\K_{s_1, t_1}^{(\lambda_1)}$, with $\gcd(s_0,t_0)=\gcd(s_1,t_1)=1$, and to simplify our analysis we assume $s_0t_0\geq1, s_1t_1\geq1,$ (identifying $\K_2$ with $\K_{1,1}$). Also we assume, for each $i$, that each component of $\Pi_i$ has $s_i$ vertices in $V_2$, so $|V_2|=r_0s_0=r_1s_1$, and $|G|=r_0s_0t_0\lambda_0=r_1s_1t_1\lambda_1$, which gives $t_0\lambda_0=t_1\lambda_1$. 

For each $i$, let $H_i$ be the subgroup of $G$ of index $r_i$. Then the $H_i$-orbits in $V_2$ are the $s_i$-subsets of vertices in the $\Pi_i$-components. Let $d=\gcd(r_0, r_1)$ and let $H$ be the index $d$ subgroup of $G$. Then  for each $i$, $H_i\leq H$ and hence the $H$-orbits in $V_2$ are unions of $H_i$-orbits, for each $i$. Let $\Delta$ be one such $H$-orbit and let $\delta\in\Delta$. Then all paths in $\Gamma$ starting from $\delta$ and ending in $V_2$ (using edges from $E_0$ or $E_1$ or both) must end at a vertex of $\Delta$. Since $\Gamma$ is connected, it follows that $\Delta=V_2$, and hence $d=1$, that is, $\gcd(r_0,r_1)=1$. Therefore $s:=s_0/r_1=s_1/r_0$ is an integer (since $r_0s_0=r_1s_1$); and each $\Pi_0$-component meets each $\Pi_1$-component in exactly $s$ vertices of $V_2$.

Let $\lambda:=\gcd(\lambda_0,\lambda_1)$ and set $t=\lambda_0/\lambda$ and $t'=\lambda_1/\lambda$, so $t_0t=t_1t'$ (since $t_0\lambda_0=t_1\lambda_1$). Since $\gcd(t,t')=1$, $u:=t_0/t'=t_1/t$ is an integer. Set $r:=r_0, r':=r_1$. Then $\Gamma=\Gamma_0^{(\lambda)}$ and each $\Pi_i=\Phi_i^{(\lambda)}$, where $\Phi_0= r\K_{sr', ut'}^{(t)}$ and $\Phi_1= r'\K_{sr, ut}^{(t')}$, and $|G|=rr'su\lambda tt'$; and we have $\gcd(r,r')=\gcd(t,t')=1$, and $\gcd(sr,ut)=\gcd(sr',ut')=1$.
Thus  line 8 of Table~\ref{t:mainbireg2} holds for $\Gamma_0$.  

This graph family covers some special cases: for example if, say,  $s=r=1$ then $\Gamma_0=\K_{r',ut'}^{(t)}+r'\K_{1,ut}^{(t')}$; and if in addition $u=t=1$ then $\Gamma_0=\K_{r',t'}+r'\K_{2}^{(t')}$.

\medskip

This completes the proof that all graphs $\Ga$ with at least three vertices, and no isolated vertices, and  admitting a cyclic edge-regular or edge-bi-regular group of automorphisms, are of the form $\Gamma=\Gamma_0^{(\lambda)}$ with $\Gamma_0$ listed in one of the Tables~\ref{t:mainreg}, \ref{t:mainbireg1}, and~\ref{t:mainbireg2}. Conversely it follows from  Lemmas~\ref{l:n-cycle}, \ref{l:circ1}, \ref{l:cyclest}, \ref{l:kst}, \ref{l:ck}, \ref{l:ck2}, \ref{l:kk}, and Corollary~\ref{c:circ1}, that all the graphs in these tables admit such groups.


\section{Proof of Theorem~\ref{thm-2}}

The aim of this section is to complete the proof of Theorem~\ref{map-cycle}. 
So let $\Ga=(V,E)$ be a graph which has a symmetrical Euler cycle 
\[
C=(e_0, e_1,\dots, e_{\ell-1})
\]
where $\ell=|E|$. Thus each $e_i=[\a_{i-1},e_i,\a_i]$ where we write $\a_\ell=\a_0$.  In particular $\Gamma$ is connected, and there exists 
$x\in\Aut\Gamma$ such that 
\[
x:\ e_i\to e_{i+2},\ \mbox{for each $i$, reading subscripts modulo $\ell$.}
\]
Then $\l x\r$ is edge-regular if $\ell$ is odd, or edge-bi-regular if $\ell$ is even, and hence, by Theorem~\ref{main-thm}, $\Gamma$ is one of the graphs in Tables~\ref{t:mainreg}, \ref{t:mainbireg1}, or~\ref{t:mainbireg2}. Our task is to decide which of these graphs has a symmetrical Euler cycle, and for each of these graphs, to determine the largest subgroup $H(C)$ of the group $D(C)$ in \eqref{e:DC} induced by $(\Aut \Gamma)_{[C]}$. Note that $x$ induces the element $\varphi^2$ of \eqref{e:varphitau} so $H(C)$ contains $\varphi^2$.
Recall that, for each vertex $\a$, the number of edges of $C$ incident with $\a$ must be even.
First we consider some of the examples listed in Theorem~\ref{map-cycle}.

\begin{lemma}\label{l:kc-euler}
\begin{enumerate}
\item[(a)]  If $\Gamma=\C_n^{(\lambda)}$ for some $n\geq3, \lambda\geq1$, then $\Gamma$ has a symmetrical Euler cycle $C$ and $H(C)=\l\varphi,\tau\r$.

\item[(b)] If $\Gamma=\K_{s,t}^{(\lambda)}$ for some $\lambda\geq 1$ and $st>1$ with $\gcd(s,t)=1$, then $\Gamma$ has a symmetrical Euler cycle $C$ if and only if $\lambda$ is even, and in this case $H(C)=\l\varphi^2,\varphi\tau\r$.
\end{enumerate}

\end{lemma}

\proof
(a)  Let $\Gamma=\C_n^{(\lambda)}$. By Proposition~\ref{p:biregext}, it is sufficient to assume that $\lambda=1$, and in this case the assertion follows from Lemma~\ref{l:n-cycle}(a).

(b) Now let $\Gamma=\K_{s,t}^{(\lambda)}$. By our comments above, for an Euler cycle to exist each vertex must be incident with an even number of edges, so both $s\lambda$ and $t\lambda$ are even. Since $\gcd(s,t)=1$, this implies that $\lambda$ is even. Therefore,
by Proposition~\ref{p:biregext}, it is sufficient to assume that $\lambda=2$ and to prove that $\Gamma=\K_{s,t}^{(2)}$ has a symmetrical Euler cycle $C$ with $H(C)$ as claimed. 
We use the notation  from Definition~\ref{d:kst} for the vertices and edges of $\K_{s,t}$ and the map $g$, and the convention in \eqref{e:Elambda} for edges of $\Gamma$. Thus the vertex set is $V=V_1\cup V_2$ with $V_1=\ZZ_s$ and $V_2=\ZZ_t$. Let $e_0=e_{0,0}^1$ and $e_1=e_{1,0}^2$, and for  each $i$ such that $1\leq i\leq st-1$, define
\[
\mbox{$e_{2i}:= e_{i,i}^1=e_0^{g^i}$ and $e_{2i+1}:= e_{i+1,i}^2=e_1^{g^i}$}.
\]
Note that $e_{2i}=[i,e_{2i},i]$ and $e_{2i+1}=[i+1,e_{2i},i]$, where we read the first entries (elements of $V_1$)  modulo $s$ and the last entries (elements of $V_2$) modulo $t$. Thus $e_{2i}, e_{2i+1}$ are both incident with $i\in V_2$, and similarly $e_{2i-1}, e_{2i}$ are both incident with $i\in V_1$. Further, since $\gcd(s,t)=1$, for each edge $e$ of $\K_{s,t}$, the edge $e^1$ occurs as $e_{2i}$ for some $i$, and $e^2$ occurs as $e_{2i+1}$ for some $i$. Also, since $|g|=st$ and $\l g\r$ is regular on the edge-set of $\K_{st}$, it follows that 
\[
C:=(e_0,e_1,\dots, e_{2st-1})
\] 
is an Euler cycle for $\Gamma$, that $C$ is preserved by $\l g\r$, and that $g$ induces the map $\varphi^2$ in $D(C)$. Thus $C$ is a symmetrical Euler cycle.   Finally, $\Aut\Gamma$ contains the following involution $y$, where $y|_{V_1}: i\to -i$ for all $i\in V_1=\ZZ_s$, $y|_{V_2}: i\to -i$ for all and all $i\in V_2=\ZZ_t$; and
\[
y: e_{i,i}^1 \to e_{-i,-i}^1,\ \mbox{and}\ y: e_{i+1,i}^2 \to e_{-i,-i-1}^2, \ \mbox{for all $i$.}
\]
It is straightforward to check that $y$ preserves $C$, namely $y:e_{2i}\to e_{2j}$ and $y:e_{2i+1}\to e_{2j-1}$ where $i+j=st$; and $y$ induces the map $\varphi\tau$ of $D(C)$ as in \eqref{e:varphitau}. Thus $H(C)$ contains  $\l\varphi^2,\varphi\tau\r$, and since $\Gamma$ is not vertex-transitive, equality holds.
\qed

Lemma~\ref{l:kc-euler} deals with all the graphs which admit a cyclic edge-regular subgroup, see Table~\ref{t:mainreg}. So we may assume that $\Gamma$ admits no cyclic edge-regular subgroup. Hence the length $\ell$ of $C$ is even, and $\l x\r$ is bi-regular on $E$, say with orbits $E_0, E_1$.  This means that, in the symmetrical Euler cycle  $C$ the edges occur alternately in $E_0$ and $E_1$. Replacing $C$ by a shift if necessary, we may assume then that $E_0=\{e_{2i}\mid 0\leq i\leq \ell-1\}$ and $E_1=\{e_{2i+1}\mid 0\leq i\leq \ell-1\}$.

\begin{lemma}\label{l:bireg-euler}
\begin{enumerate}
\item[(a)] The group $\l x\r$ has at most two orbits in $V$.
\item[(b)] If $\l x\r$ has two orbits in $V$, say $V_1$ and $V_2$, then each edge is incident with exactly one vertex of $V_1$ and one vertex of $V_2$. 
\item[(c)] Thus $\Gamma$ is one of the graphs in lines 2--4 of Table~\ref{t:mainbireg1} or line 3 of Table~\ref{t:mainbireg2}.  
\end{enumerate}
\end{lemma}

\proof
(a) Since $E_0=\{e_{2i}\mid 0\leq i\leq \ell-1\}$ and $E_1=\{e_{2i+1}\mid 0\leq i\leq \ell-1\}$, it follows from the definition of $C$ that every vertex is incident with at least one edge in $E_0$ and at least one edge in $E_1$. Then, since $\l x\r$ is transitive on $E_0$ and $E_1$ we conclude that $\l x\r$ has at most two orbits in $V$. 

(b) Suppose $\l x\r$ has two orbits $V_1$ and $V_2$ in $V$. Then as $\Gamma$ is connected, at least one of the edge-orbits, say $E_0$, has edges incident with vertices from each of $V_1$ and $V_2$, so $e_0=[\a,e_0,\b]$ with, say, $\a\in V_1$ and $\b\in V_2$. Thus $e_{2i-1}$ is incident with $\a\in V_1$ and $e_{2i+1}$ is incident with $\b\in V_2$, and hence also the edges of $E_1$ are incident with vertices from both $V_1$ and $V_2$.

(c) The last line of Table~\ref{t:mainbireg2} is not possible by part (a), graphs in line 1 of Table~\ref{t:mainbireg1} or lines 1--2 of Table~\ref{t:mainbireg2} are excluded because they admit cyclic edge-regular subgroups, and graphs in lines 4--7 of Table~\ref{t:mainbireg2} are not possible by part (b).  
\qed

The next two lemmas deal with the remaining graphs from Table~\ref{t:mainbireg1}.  

\begin{lemma}\label{Cyc-2-E}  
\begin{itemize}
\item[(a)] If $\Gamma=\Gamma_0^{(\lambda)}$ with $\Gamma_0=\Gamma(2r,1,r)=\C_{2r}+r\K_2^{(2)}$ as in Corollary~\ref{c:circ1}(a) (line $4$ of Table~\ref{t:mainbireg1}), where  $r\geq2$ and $\lambda\geq1$, then $\Gamma$ has a symmetrical Euler cycle $C$ if and only if $r$ is even, and in this case  $H(C)=\l \varphi^2, \varphi\tau\r$.

\item[(b)] If $\Gamma=\Gamma_0^{(\lambda)}$ with $\Gamma_0=\Gamma(2r,2,r)=2\C_{r}+r\K_2^{(2)}$ as in Corollary~\ref{c:circ1}(b) (line $5$ of Table~\ref{t:mainbireg1}), where $r$ is odd, $r\geq3$,  and $\lambda\geq1$, then $\Gamma$ has a symmetrical Euler cycle $C$ and  $H(C)=\l \varphi^2, \varphi\tau\r$.

\end{itemize}
\end{lemma}

\proof
We use the notation from Definition~\ref{d:circ1}.
In both cases (a) and (b), $\Gamma_0 =  \Circ(n, S)=(V,E)$ with $S=\{a, -a, r^{(2)}\}$ as in Definition~\ref{d:circ1} (with $a=1$ or $2$), and we have $V=\ZZ_{2r}$ and $4r$ edges $e_{i,a}=[i,e_{i,a},i+a]$ and $e_{i,r}=[i,e_{i,r},i+r]$, for $i\in\ZZ_{2r}$. Also the map $g:i\to i+1, e_{i,a}\to e_{i+1,a}, e_{i,r}\to e_{i+1,r}$ lies in $\Aut\Gamma_0$, by Lemma~\ref{l:circ1}, and is bi-regular on $E$ with $\l g\r$-orbits $E_a=\{e_{i,a}\mid i\in \ZZ_{2r}\}$ and $E_r=\{e_{i,r}\mid i\in \ZZ_{2r}\}$. We also use the notation from \eqref{e:Elambda} for edges of $\Gamma=(V, E^{(\lambda)})$  and from the proof of Proposition~\ref{p:biregext}(a) for the edge-bi-regular automorphism $g^{(\lambda)}$ of $\Gamma$ corresponding to the edge-bi-regular map $g$ on $\Gamma_0$. The $\l g^{(\lambda)}\r$-orbits in $E^{(\lambda)}$ are $E_a^{(\lambda)}$ and $E_r^{(\lambda)}$.

(a) Consider first case (a), so $a=1$. Then the vertex-action $(\Aut\Gamma)^V$ is contained in $(\Aut [E_1])^V\cong \D_{4r}$ and contains $g^V$ of order $2r$. 
Suppose that $C=(e_0,e_1,\dots, e_{4r\lambda-1})$ is an Euler cycle for $\Gamma$ with element $x$ inducing $\varphi^2\in H(C)$, as above. Then $\l x\r$ induces $\ZZ_{2r}$ on $V$, and hence $x^V= (g^i)^V$ for some $i$ such that $\gcd(i,2r)=1$. We show first that $r$ is even.  We may assume that the `even' edges $e_{2i}$ of $C$ lie in $E_a^{(\lambda)}$ and the `odd' edges $e_{2i+1}^{(\lambda)}$ lie in $E_r^{(\lambda)}$. Moreover, replacing $C$ by a cycle in its sequence class, if necessary, we may assume further that $e_0=e_{0,a}^1$ and $e_1$ is $e_{1,r}^u$ or $e_{r+1,r}^u$, for some $u$, and hence that $e_{2}$ is $e_{r+1,1}^v$ or $e_{r,1}^v$, for some $v$, and $e_{4r\lambda-1}$ is $e_{0,r}^w$ or $e_{r,r}^w$, for some $w$.  Considering the vertices incident with these edges, it follows from $e_0^x=e_2$  that $\{0,1\}^{g^i}=\{r+1,r+2\}$ or $\{r,r+1\}$, and from $e_{4r\lambda-1}^x=e_1$ that $\{0,r\}^{g^i}=\{1,r+1\}$. These conditions together imply that $i=r+1$, and then from $\gcd(i,2r)=1$ we conclude that $r$ is even. To complete the proof of part (a) we  exhibit a symmetrical Euler cycle $C$ for $\Gamma$ such that $H(C)=\l \varphi^2, \varphi\tau\r$ (noting that $H(C)\ne D(C)$ since $\Gamma$ is not edge-transitive).  By Proposition~\ref{p:biregext} we may assume that $\lambda=1$.    As noted above, in a symmetrical Euler cycle $C$ preserved by $x\in\Aut\Gamma$ which induces the map $\varphi^2$ in $D(C)$, we may assume that $e_0=e_{0,1}$, that $e_1=e_{1,r}$ or $e_{r+1,r}$, and that $x=g^{r+1}$. So for our construction let us choose $e_1=e_{1,r}$, so that, for all $i=1,\dots,2r-1$, 
\[
\mbox{$e_{2i}:= e_{0,1}^{g^{i(r+1)}} = e_{i(r+1),1}$ and $e_{2i+1}:=e_{1,r}^{g^{i(r+1)}}= e_{i(r+1)+1,r}$}.
\]
Then $e_{2i}, e_{2i+1}$ are both incident with $i(r+1)+1$, and $e_{2i+1}, e_{2i+2}$  are both incident with $(i+1)(r+1)$. Hence $C:=(e_0,e_1,\dots, e_{4r-1})$ is an Euler cycle of $\Gamma$. The element $g^{r+1}$ preserves the cycle $C$ and acts as $g^{r+1}: e_j\to e_{j+2}$, that is $g^{r+1}$ induces the map $\varphi^2$ of $D(C)$. The cycle $C$ is also preserved by the automorphism $y\in\Aut\Gamma$ defined by: $y: -i\leftrightarrow i+1$ on $V$ and $y: e_{-i,1}\leftrightarrow e_{i,1},$\   $e_{-i,r}\leftrightarrow e_{i+1,r}$ on $E$. This map induces $\varphi\tau\in D(C)$ as in \eqref{e:varphitau}, and hence $H(C)$ contains $\l \varphi^2, \varphi\tau\r$, and equality holds since $\Gamma$ is not edge-transitive. An exactly similar argument works for the case where $e_1=e_{r+1,r}$. This completes the proof of part (a).

(b) Now consider case (b), so $a=2$ and $r$ is odd. By Proposition~\ref{p:biregext} we may assume that $\lambda=1$. Thus, for the edge-bi-regular automorphism $g$ in the first paragraph of this proof,  $\gcd(r+2, |g|)=\gcd(r+2, 2r)=1$, and hence $x:=g^{r+2}$ generates $\l g\r$ and $\l x\r$ is edge-bi-regular on $\Gamma$.   We construct a symmetrical Euler cycle $C-(e_0,\dots, e_{4r-1})$ preserved by $x$ and on which $x$ induces $\varphi^2$ in $H(C)$. Let $e_0:= e_{0,r}$ and $e_1:=e_{r,2}$, and for $1\leq i\leq 2r-1$, let 
\[
\mbox{$e_{2i}:= e_0^{x^i}= e_{0,r}^{g^{i(r+2)}} = e_{i(r+2),r}$ and $e_{2i+1}:=e_1^{x^i} = e_{r,2}^{g^{i(r+2)}}= e_{i(r+2)+r,2}$}.
\]
Then $e_{2i}, e_{2i+1}$ are both incident with $i(r+2)+r$, and $e_{2i+1}, e_{2i+2}$  are both incident with $(i+1)(r+2)$. Hence $C$ is an Euler cycle of $\Gamma$, and by definition $x$ preserves $C$ and induces $\varphi^2$. In addition $\Aut\Gamma$ contains an involution $y$ inducing $\varphi\tau$ on $C$, namely $y:i\to r-i$ on $V=\ZZ_{2r}$, and $y:e_i\to e_{4r-i}$. This is straightforward to check. Hence $H(C)$ contains, and so is equal to $\l \varphi^2, \varphi\tau\r$.
\qed

\begin{lemma}\label{l:circ-euler}
Suppose that $\Gamma=\Gamma_0^{(\lambda)}$ with $\Gamma_0 =  \Circ(n, S)=(V,E)$ with $S=\{a, -a, b,-b\}$ as in Definition~\ref{d:circ1} with $V=\ZZ_n$, $|S|=4$ and $\gcd(n,a,b)=1$,  and that $\lambda\geq1$.
\begin{enumerate}
\item[(a)] If either  $\gcd(n,a+b)=1$ or $\gcd(n,a-b)=1$, and then $\Gamma$ has a symmetrical Euler cycle, and $H(C)=\l \varphi^2, \varphi\tau\r$ for all such cycles $C$.

\item[(b)] Moreover, if  $\Aut\Gamma_0$ is not edge-transitive, then $\Gamma$ has a symmetrical Euler cycle $C$ if and only if either $\gcd(n,a+b)=1$ or $\gcd(n,a-b)=1$.
\end{enumerate}
 \end{lemma}

\begin{remark}\label{r:circ-rem}
{\rm
We comment on part (b) of Lemma~\ref{l:circ-euler}. The proof depends on the fact that, for a symmetrical Euler cycle $C$ of $\Gamma=\Gamma_0^{(\lambda)}$ (with $\Gamma_0=\Circ(n,S)=(V,E)$ as in line 4 of Table~\ref{t:mainbireg1}), an automorphism $x$ inducing $\varphi^2\in H(C)$ acts on $V$ as a power of the map $g^V$ from Definition~\ref{d:circ1}. 
The assumption that $\Gamma_0$ is not edge-transitive is sufficient to deduce this fact. However it is possible for $\Gamma_0$ to be edge-transitive. For example if $n=5, a=1, b=2$ then $\Gamma_0\cong \K_5$. Edge-transitivity occurs more generally, for example, when $n$ is odd and there exists $d\in\ZZ_n$ with $d^2=-1$, and $b=\pm da$ so $S=\{a, da, d^2a, d^3a\}$. In this case the `multiplicative' automorphism $y:i\to di$ interchanges the two $\l g\r$ edge-orbits $E_a$ and $E_b$ defined in Definition~\ref{d:circ1}. If $\Gamma_0$ is a so-called \emph{normal Cayley graph}, that is, if $\l g, y\r = \ZZ_n\rtimes \ZZ_4$ is the full automorphism group $\Aut\Gamma_0$, then the assertion in case (b) can be proved by a slightly modified argument. However we do not know if there are parameters $n,a,b$ for which 
$\Gamma_0$ is a non-normal Cayley graph (and hence is edge-transitive) and the necessary and sufficient condition in part (b) does not hold. 
}
\end{remark} 

\proof
Here  $\Gamma_0 =  \Circ(n, S)=(V,E)$ with $S=\{a, -a, b,-b\}$ as in Definition~\ref{d:circ1} with $|S|=4$, $V=\ZZ_{n}$, and $2n$ edges $e_{i,a}=[i,e_{i,a},i+a]$ and $e_{i,b}=[i,e_{i,b},i+b]$, for $i\in\ZZ_{n}$. The map $g:i\to i+1, e_{i,a}\to e_{i+1,a}, e_{i,b}\to e_{i+1,b}$ lies in $\Aut\Gamma_0$, by Lemma~\ref{l:circ1}, and is bi-regular on $E$ with $\l g\r$-orbits $E_a=\{e_{i,a}\mid i\in \ZZ_{n}\}$ and $E_b=\{e_{i,b}\mid i\in \ZZ_{n}\}$. We also use the notation from \eqref{e:Elambda} for edges of $\Gamma=(V, E^{(\lambda)})$. 

(a) Suppose first that $\gcd(a+b,n)=1$ or $\gcd(a-b,n)=1$. We construct a symmetrical Euler cycle for $\Gamma$. By Proposition~\ref{p:biregext} it is sufficient to do this for $\lambda=1$, that is, for $\Gamma=\Gamma_0=\Circ(n,S)$. So let $j=a\pm b$ and assume that $\gcd(j,n)=1$. First let $j=a+b$, let $e_0=e_{0,a}, e_1=e_{a,b}$ and for $1\leq i\leq n-1$, let
\[
\mbox{$e_{2i}:= e_0^{g^{ij}}=  e_{ij,a}$ and $e_{2i+1}:=e_1^{g^{ij}} =  e_{ij+a,b}$}.
\]
Then $C=(e_0,\dots,e_{2n-1})$ is a symmetrical Euler cycle and $g^j$ acts as $\varphi^2\in H(C)$. In addition $y\in\Gamma$ defined by $y:i\to a-i$ on $V$, and $y:e_i\to e_{2n-i}$, induces $\varphi\tau$ on $C$ so $H(C)\geq \l \varphi^2,\varphi\tau\r$. 
Similarly, for $j=a-b$, let $e_0=e_{0,a}, e_1=e_{a-b,b}$ (with $e_1$ incident with $\{a,a-b\}$) and for $1\leq i\leq n-1$, let
\[
\mbox{$e_{2i}:= e_0^{g^{ij}}=  e_{ij,a}$ and $e_{2i+1}:=e_1^{g^{ij}} =  e_{ij,b}$}.
\]
Then again $C=(e_0,\dots,e_{2n-1})$ is a symmetrical Euler cycle and $g^j$ acts as $\varphi^2\in H(C)$.  Again the element $y\in\Gamma$ defined by $y:i\to a-i$ on $V$ and $y:e_i\to e_{2n-i}$ induces $\varphi\tau$ on $C$ so $H(C)\geq\l \varphi^2,\varphi\tau\r$.
If, for some symmetrical Euler cycle $C$, $H(C)$ is strictly larger than $\l \varphi^2,\varphi\tau\r$, then $H(C)=D(C)$ and so some element $x\in \Aut\Gamma$ induces $\varphi$ on $C$, and hence $\l x\r$ is a cyclic edge-regular subgroup of $\Aut\Gamma$. This however  contradicts Theorem~\ref{main-thm} as $\Gamma=\Gamma_0$ does not appear in Table~\ref{t:mainreg}.  Thus part (a) is proved. 

(b) Suppose now that $\Aut\Gamma_0$ is not edge-transitive and that $\Gamma$ has a symmetrical Euler cycle $C=(e_0,\dots,e_{2n\lambda-1})$ and $x\in\Aut\Gamma$ preserves $C$ and acts by $x:e_i\to e_{i+2}$ for each $i$. Since $\Gamma_0$ is not edge-transitive, $\Aut\Gamma_0\leq \Aut [E_a]$ and hence the groups induced on $V$ by $\Aut\Gamma$ and $\Aut\Gamma_0$ are both isomorphic to $\D_{2n}$ with index $2$ subgroup $\l g^V\r$ generated by $g^V:i\to i+1$, with $g$ as in the first paragraph of the proof. Since $\l x\r$ is cyclic of order $n\lambda$ it follows that $x^V$ has order $n$ and hence $x^V = (g^V)^j$ for some $j$ coprime to $n$.   Relabelling the elements of $S$ if necessary, we may assume that $e_0=e_{0,a}^u$ for some $u$. Then, the next edge $e_1$ lies in $E_b^{(\lambda)}$ and hence is $e_{a,b}^v$ or $e_{a-b,b}^v$ for some $v$. Assume first that $e_1=e_{a,b}^v$. Then $e_2$ is incident with $a+b$ while, since $e_2=e_0^x$, the pair of vertices incident with $e_2$ are $\{0,a\}^x=\{0,a\}^{g^j}=\{j,a+j\}$. Hence $j=a+b$ or $j=b$ and $e_2=e_{a+b,a}^w$ or $e_{b,a}^w$, for some $w$, and $e_3$ is incident with $2a+b$ or $b$, respectively. Moreover since $e_3=e_1^x$, the pair of vertices incident with $e_3$ are $\{a,a+b\}^x=\{a,a+b\}^{g^j}$ which is $\{2a+b,2a+2b\}$ if $j=a+b$ or $\{a+b,a+2b\}$ if $j=b$. However, if $j=b$ then $e_3$ should be incident with $b$, and $b\not\in\{a+b,a+2b\}$ since $a\ne 0,-b$. Hence $j=a+b$ so $\gcd(a+b,n)=1$. A similar argument in the case where $e_1=e_{a-b,b}^v$ leads to the condition $\gcd(a-b,n)=1$. This proves the asserted congruence conditions, and the converse follows from part (a).
This completes the proof.
\qed

Finally we consider the remaining graphs from Table~\ref{t:mainbireg2}.  

\begin{lemma}\label{l:kst-euler}
If $\Gamma=\Gamma_0^{(\lambda)}$ with $\Gamma_0=\C_{2r}[s\K_1,t\K_1]$ as in Definition~\ref{d:cyclest} (line $4$ of Table~\ref{t:mainbireg2}),  with $r\geq2, st\geq2$, $\gcd(s,t)=1$, and $\lambda\geq1$, then $\Gamma$ has a symmetrical Euler cycle $C$ and  $H(C)=\l \varphi^2, \tau\r$.
\end{lemma}
 
\proof
We use the notation from Definition~\ref{d:cyclest}.
Thus $\Gamma_0 = \C_{2r}[s\K_1,t\K_1]=(V,E)$, and $V=V_S\cup V_T$ and $E=E_S\cup E_T$, where $S=\ZZ_s, T=\ZZ_t$. Also we have $g, y\in\Aut \Gamma_0$ such that $\l g\r$ is edge-bi-regular with edge-orbits $E_S, E_T$ and $\l g,y\r\cong \D_{2rst}$ is bi-transitive on $V$ with vertex-orbits $V_S, V_T$, and regular on $E$ (see Lemma~\ref{l:cyclest}).
By Proposition~\ref{p:biregext} it is sufficient to construct a cycle $C$ with the required properties in the case $\lambda=1$. So we assume now that $\lambda=1$ and $\Gamma=\Gamma_0$, and we construct $C=(e_0,\dots, e_{2rst-1})$. 

From Definition~\ref{d:cyclest} the edges in $E_S, E_T$ are labelled $e^{2k}_{i,j}$ and $e^{2k+1}_{j,i}$, respectively, where $0\leq k\leq r-1, i\in S, j\in T$. Each integer $i$ such that $0\leq i\leq 2rst-1$ can be represented uniquely as $i=\ell+2rm$, where $0\leq \ell\leq 2r-1$ and $0\leq m\leq st-1$. For such a representation we define $e_i=e^\ell_{m,m}$ if $\ell\leq 2r-2$ and  $e_i=e^\ell_{m,m+1}$ if $\ell=2r-1$. If $\ell\leq 2r-2$, then $e_i, e_{i+1}$ are both incident with the vertex $(\ell,m)$ where $m\in T$ if $\ell$ is even or $m\in S$ if $\ell$ is odd. If $\ell= 2r-1$, then $e_i, e_{i+1}$ are both incident with the vertex $(0,m+1)$ with $m+1\in S$. This proves that $C$ is a cycle of length $2rst$, and from its definition $C$ involves all edges and hence is an Euler cycle. Moreover $C$ is preserved by $g$ and $y$, with $g$ inducing $\varphi^2\in H(C)$ and $y$ inducing  $\tau\in H(C)$. If $H(C)$ is strictly larger than $\l \varphi^2, \tau\r$, then $H(C)=D(C)$, which implies that $\Gamma_0$ is vertex-transitive, and hence that $s=t$, which is a contradiction. This completes the proof.
\qed


\begin{thebibliography}{xx}

\bibitem{DX}
S. Du and M-Y Xu, A classification of semisymmetric graphs of order $2pq$, {\it Comm. Alg.} 
{\bf 28} (2000), 2685--2715. 

\bibitem{1-face}
W. W. Fan, C. H. Li and Naer Wang,
Edge-transitive uniface embeddings of complete bipartite multi-graphs,
{\it J. Algebraic Combin.} {\bf 49} (2019), 125--134. 

\bibitem{GLP}
M. Giudici, C. H. Li and C. E. Praeger, Analysing finite locally $s$-arc transitive graphs,  
\emph{Trans. Amer. Math. Soc.} {\bf 356} (2004), 291--317.



\bibitem{LPS-rotary}
C. H. Li, C. E. Praeger and S-J Song, Vertex-rotary maps and coset graphs with finite edge multiplicity, preprint 2021.

\bibitem{BigMaps}
C. H. Li, C. E. Praeger and S-J Song, Arc-transitive embeddings of graphs, preprint 2021.

\bibitem{ci-map}
R. B. Richter, P.D. Seymour and J. \v Sir\'a\v n, Circular embeddings of planar graphs in non-spherical
surfaces. {\it Discrete Math}. {\bf 126} (1994), 273--280.

\bibitem{C-map}
N. Robertson, X. Y. Zha, Closed 2-cell embeddings of graphs with no V8-minors,  {\it Discrete Math}.
{\bf 230} (2001), 207--213.

\bibitem{Singerman}
D. Singerman, Unicellular dessins and a uniqueness theorem for Klein's Riemann surface of genus 3,
{\it Bull. London Math. Soc.} {\bf 33} (2001), 701--710.

\bibitem{s-map}
X. Y. Zha, The closed 2-cell embeddings of 2-connected doubly toroidal graphs.
{\it Discrete Math}. {\bf 145} (1995), 259--271.
\end{thebibliography}
\end{document}